\documentclass[final]{siamltex}
\usepackage{epsfig,amsmath,amsfonts,amssymb,graphicx,euscript}
\usepackage{palatino}

\newtheorem{remark}[theorem]{Remark}
\newtheorem{example}[theorem]{Example}
\newcommand{\real}{\mathbb{R}}

\newcommand{\mathscr}[1]{\boldsymbol{\EuScript{#1}} }

\newcommand{\spn}{\operatorname{span}}
\newcommand{\rk}{\operatorname{rk}}
\newcommand{\id}{\operatorname{Id}}

\newcommand{\Hc}{{\cal H}}
\newcommand{\Rc}{{\cal R}}

\newcommand{\Sc}{\mathcal{S}}
\newcommand{\Uc}{\mathcal{U}}
\newcommand{\Nc}{\mathcal{N}}
\newcommand{\Yc}{\mathcal{Y}}

\newcommand{\Lc}{\mathcal{L}}
\newcommand{\Zc}{\mathcal{Z}}
\newcommand{\D}{{\mathcal{D}}}
\newcommand{\G}{{\mathcal{G}}}
\newcommand{\eps}{\epsilon}

\def\frak{\mathfrak{X}}

\newcommand{\setdef}[2]{\{#1 \; | \; #2\}}

\newcommand{\supscr}[2]{#1^{\textsc{#2}}}
\newcommand{\symprod}[2]{\langle #1: #2\rangle}
\newcommand{\Beltrami}[2]{\{ #1: #2 \}}

\newcommand{\pder}[2]{\frac{\partial #1}{\partial #2}}

\newcommand\map[3]{#1\colon#2\rightarrow#3}

\renewcommand{\theenumi}{(\roman{enumi})}

\title{Characterization of gradient control systems\footnotemark[1]}

\author{Jorge Cort\'es\footnotemark[2]
\and Arjan van der Schaft\footnotemark[3]
\and Peter E. Crouch\footnotemark[4]
}

\begin{document}
\maketitle 

\renewcommand{\thefootnote}{\fnsymbol{footnote}}

\footnotetext[1]{A short version of this paper was presented as~\cite{CoScCr}
  in the IFAC Workshop on Lagrangian and Hamiltonian Methods for Nonlinear
  Control, Seville, Spain, 2003.}

\footnotetext[2]{Coordinated Science Laboratory, University of Illinois at
  Urbana-Champaign, 1308 W. Main St., Urbana, IL 61801, United States, Ph.
  +1 217 244-8734, Fax. +1 217 244-1653, \texttt{jcortes@uiuc.edu},
  \texttt{http://motion.csl.uiuc.edu/\~{}jorge}}

\footnotetext[3]{Department of Applied Mathematics, University of Twente, PO
  Box 217, 7500 AE Enschede, The Netherlands, Ph. +31 53 489-3449, Fax +31 53
  434-0733, \texttt{a.j.vanderschaft@math.utwente.nl},
  \texttt{http://www.math.utwente.nl/\~{}twarjan}}

\footnotetext[4]{Department of Electrical and Computer Engineering, Arizona
  State University, Tempe, AZ 85287, United States, Ph. +1 480 965-1722, Fax.
  +1 480 965-2267, \texttt{peter.crouch@asu.edu},
  \texttt{http://www.eas.asu.edu/\~{}sserc/people/crouch/crouch.html}}

\renewcommand{\thefootnote}{\arabic{footnote}}

\begin{abstract}
  We investigate necessary and sufficient conditions under which a general
  nonlinear affine control system with outputs can be written as a gradient
  control system corresponding to some pseudo-Riemannian metric defined on the
  state space.  The results rely on a suitable notion of compatibility of the
  system with respect to a given affine connection, and on the output behavior
  of the prolonged system and the gradient extension.  The symmetric product
  associated with an affine connection plays a key role in the discussion.
\end{abstract}



\pagestyle{myheadings} \thispagestyle{plain} \markboth{J. Cort\'es, A.J. van
  der Schaft, P.E. Crouch}{Characterization of gradient control systems}

\section{Introduction}

A physically motivated class of nonlinear systems are \emph{gradient control
  systems}, see~\cite{BrMo,Cr,Sc,Sc2,Va,Ve,Wi} and the references quoted
therein. These systems are described in the following way: they are nonlinear
affine control systems, which are endowed with a pseudo-Riemannian metric on
the state space manifold. The drift vector field of the system is the gradient
vector field associated with an internal potential function with respect to
the pseudo-Riemannian metric, and the input vector fields are the gradient
vector fields associated with the output functions of the system.  Examples of
gradient control systems include nonlinear electrical RLC networks, and
dissipative systems where the inertial effects are neglected.  In the case of
RL or RC networks, the pseudo-Riemannian metric is positive-definite, and thus
is a usual Riemannian metric, while for general RLC networks the metric is
indefinite.  We refer to~\cite{BrMo,Cr,Va,Ve} for more background on the
modeling of nonlinear networks as gradient systems.

Another relevant class of nonlinear systems is the family formed by
the \emph{Hamiltonian control systems}. In this case, the state space
manifold is equipped with a symplectic form. The drift vector field
and the input vector fields are the Hamiltonian vector fields
associated, respectively, to an internal energy function and the
output functions of the system with respect to the symplectic form.
Hamiltonian equations are of central importance in the modeling of
physical systems as they are the starting point to describe the
dynamics of a very large class of phenomena, including mechanical,
electrical and electromagnetic systems.

Apart from their physical and engineering importance, gradient and
Hamiltonian systems also possess very peculiar mathematical
properties. For instance, a linear input-state-output system is a
Hamiltonian control system~\cite{BrRa} (respectively a gradient
control system~\cite{Wi}) if and only if its impulse response matrix
$W(t)$ satisfies $W(t) = -W^T(-t)$ (respectively $W(t) = W^T(t)$).
Although they are typically not amenable to linearization techniques,
their rich geometric structure makes possible to combine powerful
tools from Nonlinear Control Theory, Differential Geometry and
Classical Mechanics in the study of a variety of problems including
stability and stabilization, input-output decoupling, structural
synthesis and interconnection.

Their theoretical and practical relevance, together with their
meaningful geometric properties and the wide range of results
available for them, make the classes of Hamiltonian and gradient
systems \emph{distinct} within the family of nonlinear affine control
systems. This explains the interest in identifying those systems that
can be written as either Hamiltonian or gradient. This
\emph{characterization problem} is motivated by the Realization
Problem in Systems Theory and the Inverse Problem in Mechanics.  The
Realization Problem addresses the question of when the input-output
map of a system can be realized as the external behavior of a
Hamiltonian (respectively gradient) input-output system.  The Inverse
Problem, which has a longstanding history in mathematical physics,
poses the question of when a second-order differential equation can be
realized as the Euler-Lagrange equations corresponding to certain
Lagrangian function. For further reference on these problems, the
reader is referred to~\cite{CrIr,Ja,Sa,SaThPr}.

In~\cite{CrLaSc,CrSc}, necessary and sufficient conditions were given under
which a minimal nonlinear affine control system with an equal number of inputs
and outputs is a Hamiltonian control system with respect to some symplectic
structure, which turned out to be unique. In the present paper we describe an
analogous theory for the case of gradient systems.  As we discuss below, there
are a number of key differences in the treatment of the characterization
problem for the Hamiltonian and the gradient case, which make the latter more
involved. The role played in the Hamiltonian setting by the Lie bracket and
the Hamiltonian vector fields is taken in the gradient setting by the
symmetric product associated with an affine connection and the gradient vector
fields. A fundamental observation is that, while every input-state-output
system admits a natural extension to a Hamiltonian system living on the
cotangent bundle of its state space, the construction of a gradient extension
on the cotangent bundle relies on the selection of a torsion-free connection
on the state space. This motivates the introduction of a novel
\emph{compatibility condition} of the given nonlinear system with the selected
affine connection guaranteeing an appropriate choice of the latter one. The
compatibility condition is expressed as a relation of the symmetric products
of the drift vector field and the input vector fields with the output
functions of the system, and turns out to play a prominent role in the
characterization of gradient control systems presented in
Theorem~\ref{th:characterization} below.

The paper is organized as follows.  In
Section~\ref{se:gradient-control-systems} we present the class of nonlinear
systems considered along the paper.  We also introduce the notions of
prolongation and gradient extension of a nonlinear system, whose observability
properties are studied in Section~\ref{se:observability}.
Section~\ref{se:externally-equivalent} introduces the concept of (weakly)
externally equivalent systems.  In Section~\ref{se:main}, we introduce the
important notion of compatibility between a nonlinear system and a given
affine connection.  At this point, we are ready to state and prove the main
result of the paper, namely the characterization of when a general nonlinear
control system is gradient.  In Section~\ref{se:uniqueness-realization} we
investigate the uniqueness (up to isometry) of gradient realizations with the
same input-output behavior and we give an alternative proof of a result
in~\cite{Ba,Ba2}. We present our conclusions in Section~\ref{se:conclusions}.
Finally, an appendix in Section~\ref{se:appendix} contains a simplifying
result concerning the checkability of the compatibility condition for a
nonlinear affine control system.

\section{Setting}\label{se:gradient-control-systems}

Let $M$ be an $n$-dimensional differentiable manifold. We will denote by $TM$,
$T^*M$ the tangent and cotangent bundles of $M$, by $\mathfrak{X} (M)$ the set
of smooth vector fields on $M$, by $\Omega^1 (M)$ the set of smooth one-forms
on $M$, and by $C^{\infty}(M)$ the set of smooth functions on $M$.  Throughout
the paper, the manifold $M$ and the mathematical objects defined on it will be
assumed to be real-analytic.

Consider a nonlinear control system $\Sigma$ with state space $M$, affine in
the inputs, and with an equal number of inputs and outputs,
\begin{equation}\label{eq:system}
  \Sigma : \left\{
    \begin{array}{l}
      \displaystyle{\dot{x} = g_0 (x) + \sum_{j=1}^m u_j g_j (x)} \, , \\
      y_j = V_j(x) \, , \quad j=1,\ldots,m \, , 
    \end{array}
  \right.
\end{equation}
where $x \in M$, $x(0) = x_0$ and $u = (u_1,\ldots,u_m) \in U \subset
\real^m$. The vector fields $g_0,g_1,\ldots,g_m$ on $M$ are assumed to be
complete and $V_1,\ldots,V_m$ are real-valued functions on $M$. The set $U$ is
the control space, which for simplicity is assumed to be an open subset of
$\real^m$, containing $0$. The function $ u(t) = ( u_1(t), \ldots, u_m(t))$
belongs to a certain class of functions of time, denoted by $\Uc$, called the
\emph{admissible controls}.  For our purposes, we may restrict the admissible
controls to be the piecewise constant right continuous functions.

An important subclass of the family of nonlinear systems~\eqref{eq:system} is
formed by the Hamiltonian control systems, see~\cite{CrSc}. Here, we will
instead focus our attention on the family of gradient control systems. Let
$\G$ be a pseudo-Riemannian metric on $M$, i.e. a non-degenerate symmetric
(0,2)-tensor on $M$ (not necessarily positive definite). Consider the
`musical' isomorphisms associated with $\G$, $\flat_\G : \mathfrak{X} (M)
\rightarrow \Omega^1 (M)$, $\sharp_\G : \Omega^1 (M) \rightarrow \mathfrak{X}
(M)$ defined by
\begin{align*}
  \flat_\G (X) (Y) = \G (X,Y) \, , \quad \sharp_\G (\omega) = \flat_\G^{-1}
  (\omega) \, ,
\end{align*}
where $X, Y \in \mathfrak{X}(M)$ and $\omega \in \Omega^1(M)$. The
\emph{gradient vector field} associated with a function $V \in C^{\infty}(M)$
is given by $\grad_\G V = \sharp_\G (dV)$. Reciprocally, a vector field $X
\in \frak (M)$ is said to be \emph{locally gradient} if the one-form
$\flat_\G (X)$ is closed. By Poincar\'e's lemma, this is equivalent to saying
that there exists a locally defined function $V \in C^{\infty}(M)$ such that
$\flat_\G (X) = dV$. If this equality holds globally, $X$ is called
\emph{gradient} and will be denoted by $X = \grad_\G V$.  Along the paper, we
will drop the subindex when it is clear from the context the
pseudo-Riemannian metric with respect to which the gradient vector field is
computed. If we fix coordinates $(x^1,\ldots,x^n)$ on $M$, then the
pseudo-Riemannian metric can be locally expressed as $\G = \G_{ab} dx^a
\otimes dx^b$, where $(\G^{ab} = \G (\pder{}{x^a},\pder{}{x^b}))$ is a
symmetric matrix. The musical isomorphisms are then given by $\flat_\G =
\G_{ab} dx^a \otimes dx^b$, $\sharp_\G = \G^{ab} \pder{}{x^a} \otimes
\pder{}{x^b}$, where $(\G^{ab})$ is the inverse matrix of $(\G_{ab})$.
Finally, the gradient vector field associated with $V$ reads
\[
\grad_\G V = \G^{ab} \pder{V}{x^b} \pder{}{x^a} \, .
\]

Now, assume that the state space $M$ in~\eqref{eq:system} is a
pseudo-Riemannian manifold, $(M, \G)$. Furthermore, assume that the drift
vector field $g_0$ is locally gradient and the input vector fields $g_j$,
$j=1,\ldots,m$ are gradient with respect to the functions $V_1, \ldots, V_m$,
i.e.  $g_j = \grad_\G V_j$, $j =1,\ldots,m$. Then, the resulting system
\begin{equation}\label{eq:gradient-system}
  \Sigma : \left\{
    \begin{array}{l}
      \displaystyle{\dot{x} = g_0 (x) + \sum_{j=1}^m u_j(t) \grad_\G
      V_j (x)} \, , 
      \\ 
      y_j = V_j(x) \, , \quad j=1,\ldots,m \, ,
    \end{array}
  \right.
\end{equation}
is called a \emph{locally gradient control system on $M$}. If the drift $g_0$
is a gradient vector field, then the system is called a \emph{gradient
  control system on $M$}.

Our objective is to characterize when a nonlinear system of the
form~\eqref{eq:system} is actually a locally gradient control
system~\eqref{eq:gradient-system}, i.e. find necessary and sufficient
conditions for the existence of a pseudo-Riemannian metric on the state space
$M$ such that the system~\eqref{eq:system} equals
system~\eqref{eq:gradient-system}. These conditions will be given in terms of
the output behavior of the so-called prolonged system and the gradient
extension of $\Sigma$, which we describe next.

\subsection{The prolongation of a nonlinear system}

Given an initial state $x(0)=x_0$, take a coordinate neighborhood of $M$
containing $x_0$. Let $x(t)$, $t \in [0,T]$ be the solution
of~\eqref{eq:system} corresponding to the input functions
$u(t)=(u_1(t),\ldots,u_m(t))$ and the initial state $x(0)=x_0$, such that
$x(t)$ remains within the selected coordinate neighborhood.  Denote the
resulting output by $y(t)=(y_1(t),\ldots,y_m(t))$, with $y_j(t)=V_j(x(t))$.
Then the \emph{variational system} along the state-input-output trajectory
$(x(t),u(t),y(t))$ is given by the following time-varying system,
\begin{align}\label{eq:variational}
  & \dot{v}(t) = \pder{g_0}{x} (x(t)) v(t) + \sum_{j=1}^m u_j(t)
  \pder{g_j}{x}(x(t)) v(t) + \sum_{j=1}^m u_j^p g_j(x(t)) \, , \nonumber
  \\
  & \supscr{y_j}{v} (t) = \pder{V_j}{x} (x(t)) v(t) \, , \quad j=1,\ldots,m
  \, ,
\end{align}
where $v(0) = v_0 \in \real^n$, and $u^p = (u_1^p, \ldots, u_m^p)$,
$\supscr{y}{v} = (\supscr{y_1}{v} ,\ldots, \supscr{y_m}{v} )$ denote the
inputs and the outputs of the variational system. The reason behind the
terminology `variational' comes from the following fact: let
$(x(t,\eps),u(t,\eps),y(t,\eps))$, $t\in [a,b]$ be a family of
state-input-output trajectories of~\eqref{eq:system} parameterized by $\eps
\in (-\delta, \delta)$, with $x(t,0)=x(t)$, $u(t,0)=u(t)$ and $y(t,0)=y(t)$,
$t \in [a,b]$. Then, the infinitesimal variations
\[
v(t) = \pder{x}{\eps}(t,0) \, , \quad u^p (t) = \pder{u}{\eps}(t,0) \, ,
\quad \supscr{y}{v}(t) = \pder{y}{\eps}(t,0) \, ,
\]
satisfy equation~\eqref{eq:variational}. Additionally, if the initial state is
the same for the whole family of trajectories, $x(0,\eps)=x_0$, then the
variational state $v(0)$ at time $0$ is necessarily $0$.

The \emph{prolongation or prolonged system} of~\eqref{eq:system} corresponds
to considering together the original system~\eqref{eq:system} and the
variational system,
\begin{align}\label{eq:prolonged-system-coordinates}
  & \dot{x} = g_0 (x) + \sum_{j=1}^m u_j g_j (x) \, , \nonumber \\
  & \dot{v}(t) = \pder{g_0}{x} (x(t)) v(t) + \sum_{j=1}^m u_j(t)
  \pder{g_j}{x}(x(t)) v(t) + \sum_{j=1}^m u_j^p g_j(x(t)) \, , \nonumber
  \\
  & y_j = V_j(x) \, , \quad \supscr{y_j}{v} (t) = \pder{V_j}{x} (x(t)) \,
  v(t) \, , \quad j=1,\ldots,m \, ,
\end{align}
with inputs $u_j$, $u_j^p$, outputs $y_j$, $\supscr{y_j}{v}$ and state
$(x,v)$. To state a coordinate-free definition of the prolonged
system~\eqref{eq:prolonged-system-coordinates} on the whole tangent space
$TM$, we need to introduce the notions of vertical and complete lifts of
functions and vector fields. We do this following~\cite{YaIs}. Given a
function $V$ on $M$, the \emph{complete lift} of $V$ to $TM$,
$\supscr{V}{c}:TM \rightarrow \real$ is defined by $\supscr{V}{c} (v) =
\langle dV,v \rangle$. In the induced local coordinates on $TM$,
$(x^1,\ldots,x^n, v^1,\ldots,v^n)$, this reads
\[
\supscr{V}{c} (x,v) = \sum_{a=1}^n \pder{V}{x^a}(x) \, v_a \, .
\]
The \emph{vertical lift} of $V$ to $TM$, $\supscr{V}{v} : TM \rightarrow
\real$, is defined by $\supscr{V}{v} = V \circ \tau_M$, where $\tau_M$ denotes
the tangent bundle projection. Given a vector field $X$ on $M$, the
\emph{complete lift} of $X$ to $TM$, $\supscr{X}{c} \in \frak (TM)$ is defined
as the unique vector field verifying $\supscr{X}{c} (\supscr{f}{c}) =
\supscr{(Xf)}{c}$, for any $f \in C^{\infty}(M)$.  Alternatively, if $\Phi_t
:M \rightarrow M$, $t \in [0,\eps)$ denotes the flow of $X$, then we can
define $\supscr{X}{c}$ as the vector field whose flow is given by $(\Phi_t)_*:
TM \rightarrow TM$. In local coordinates,
\begin{gather}\label{eq:complete-lift}
  \supscr{X}{c} (x,v) = \sum_{a=1}^n X_a(x) \pder{}{x^a} + \sum_{a,b=1}^n
  \pder{X_a}{x^b} (x) v^b \pder{}{v^a} \, .
\end{gather}
The \emph{vertical lift} of $X$ to $TM$, $\supscr{X}{v} \in \frak (TM)$ is the
unique vector field such that $\supscr{X}{v} (\supscr{f}{c}) =
\supscr{(Xf)}{v} $, for any $f \in C^{\infty}(M)$. In local coordinates,
\begin{gather}\label{eq:vertical-lift}
  \supscr{X}{v} (x,v) = \sum_{a=1}^n X_a(x) \pder{}{v^a} \, .
\end{gather}

The following definition provides an intrinsic way of pasting together the
system~\eqref{eq:system} with the variational systems associated with its
state-input-output trajectories.

\begin{definition}\label{dfn:prolonged-system}
  The prolonged system $\Sigma^p$ of a nonlinear system $\Sigma$ of the
  form~\eqref{eq:system} is defined by
  \begin{equation}\label{eq:prolonged-system}
    \Sigma^p : \left\{
      \begin{array}{l}
        \displaystyle{\dot{x}_p = \supscr{g_0}{c} (x_p) + \sum_{j=1}^m u_j(t)
        \supscr{g_j}{c} (x_p) + \sum_{j=1}^m u_j^p(t)
        \supscr{g_j}{v} (x_p)} \\ 
        y_j = \supscr{V_j}{v} (x_p) \, , \quad 
        \supscr{y_j}{v} =\supscr{V_j}{c} (x_p) \, , \quad j=1,\ldots,m \, ,
      \end{array}
    \right.
  \end{equation}
  where $x_p = (x,v) \in TM$, and $x_p(0) = (x_0,v_0)$.
\end{definition}

One can easily check that in the induced tangent bundle coordinates, the local
expression of the system~\eqref{eq:prolonged-system} is
precisely~\eqref{eq:prolonged-system-coordinates}.

\begin{remark}
  {\rm In the same way as we have presented above, one can also introduce the
    notions of adjoint variational system and Hamiltonian extension of the
    nonlinear system~\eqref{eq:system}. These notions play a key role in the
    characterization of when a general system admits a Hamiltonian
    description, see~\cite{CrSc}.}
\end{remark}

\subsection{The gradient extension of a nonlinear system}

When dealing with the Hamiltonian extension of a nonlinear system, one relies
on the fact that the cotangent bundle is endowed with a canonical symplectic
structure. However, this is not the case when treating gradient systems,
since a canonical pseudo-Riemannian structure on the cotangent bundle does
not exist. In order to define the gradient extension of a nonlinear system of
the form~\eqref{eq:system}, we will first select a torsion-free affine
connection $\nabla$ on $M$, and then consider its Riemannian extension to
$T^*M$ (cf.~\cite{PaWa}).

Let us briefly present some basic notions on affine connections and
Riemannian geometry. An \emph{affine connection}~\cite{KoNo} on a manifold
$M$ is defined as an assignment
\[
\setlength{\arraycolsep}{2pt}
\begin{array}{rccc}
  \nabla: & \mathfrak{X}(M) \times \mathfrak{X} (M) & \longrightarrow & 
  \mathfrak{X} (M) \\
  & (X , Y) & \longmapsto & \nabla_X Y
\end{array}
\]
which is $\real$-bilinear and satisfies $\nabla_{fX} Y =f \nabla_X Y$ and
$\nabla_X(fY)=f \nabla_XY + X(f)Y$, for any $X$, $Y \in \mathfrak{X}(M)$,
$f\in C^{\infty}(M)$. This implies that $\nabla_X Y (x)$ only depends on
$X(x)$ and the value of $Y$ along a curve which is tangent to $X$ at $x$. Let
$c:t \in [t_0,t_1] \mapsto c(t)= (x^1(t),\ldots,x^n(t))\in M$ be a curve on
$M$ and $W$ a vector field along $c$, i.e. a map $\map{W}{[t_0,t_1]}{TM}$
such that $\tau_M(W(t))=c(t)$ for all $t\in[a,b]$. Let $V$ be a vector field
that satisfies $V(c(t))=W(t)$.  The \emph{covariant derivative of $W$ along
  $c$} is defined by
\begin{equation*}
  \frac{DW(t)}{dt}=\nabla_{\dot{c}(t)}W(t)=
  \nabla_{\dot{c}(t)}V(x)\big|_{x=c(t)}. 
\end{equation*}
This definition makes sense because of the defining properties of the affine
connection. Now, we may take $W(t)=\dot{c}(t)$ and set up
$\nabla_{\dot{c}(t)} \dot{c}(t) = 0$. This equation is called the
\emph{geodesic equation}, and its solutions are termed the \emph{geodesics}
of~$\nabla$. In local coordinates, this condition can be expressed as
$\ddot{x}^a + \Gamma^a_{bc}(x) \dot{x}^b \dot{x}^c =0$, $1 \le a \le n $,
where the $\Gamma^a_{bc}(x)$ are the \emph{Christoffel symbols} of the affine
connection, defined by
\[
\nabla_{\frac{\partial}{\partial x^b}} \frac{\partial}{\partial x^c} =
\Gamma^a_{bc} (x) \frac{\partial}{\partial x^a} \, .
\]
The vector field $S$ on $TM$ describing the geodesic equation is called the
\emph{geodesic spray} associated with the affine connection $\nabla$. In
local coordinates,
\[
S=v^a \frac{\partial}{\partial x^a} - \Gamma^a_{bc} (x) v^b
v^c\frac{\partial}{\partial v^a} \, .
\]
Therefore, the integral curves of the geodesic spray $S$ are the solutions of
the geodesic equation. The torsion tensor of an affine connection is defined
by
\[
\setlength{\arraycolsep}{2pt}
\begin{array}{rccl}
  T: & \frak (M) \times \frak (M) & \longrightarrow & \frak (M) \\
  & (X,Y) & \longmapsto & \nabla_X Y - \nabla_Y X - [X,Y] \, ,
\end{array}
\]
Locally, we have
\[
T(\pder{}{x^a},\pder{}{x^b}) = \left( \Gamma_{ab}^c - \Gamma_{ba}^c \right)
\pder{}{x^c} \, ,  
\]
Given an affine connection, the \emph{symmetric product}~\cite{LeMu} of two
vector fields $X, Y \in \frak (M)$ is defined by the operation
\[
\left \langle X:Y\right \rangle = \nabla_XY + \nabla_YX \, .
\]
The geometric meaning of the symmetric product is the following~\cite{Le}: a
distribution $\D$ on $M$ is geodesically invariant (meaning that each geodesic
of $\nabla$ whose initial velocity is in $\D$ has all its velocities in $\D$)
if and only if $\symprod{X}{Y} \in {\D}$, for all $X$, $Y \in {\D}$. The
symmetric product plays a crucial role within the so-called affine connection
formalism for mechanical control systems in the study of a variety of aspects
such as controllability, series expansions, motion planning and optimal
control~\cite{BuLe}.

Associated with the metric $\G$ there is a natural affine connection, called
the \emph{Levi-Civita} connection. The Levi-Civita connection $\nabla^\G$ is
determined by the formula
\begin{multline*}
  2 \, \G(\nabla^\G_XY,Z) = X(\G(Y,Z)) +Y(\G(Z,X)) -Z(\G(X,Y)) \\
  + \G(Y,[Z,X]) - \G(X,[Y,Z]) + \G(Z,[X,Y]) \, , \quad X,Y,Z \in
  \mathfrak{X}(M) \, .
\end{multline*}
One can compute the Christoffel symbols of $\nabla^{\G}$ to be
\begin{gather}\label{eq:Christoffel}
  \Gamma^a_{bc} = \frac{1}{2} \G^{ad} \left( \frac{\partial \G_{db}}{\partial
      x^c} + \frac{\partial \G_{dc}}{\partial x^b}- \frac{\partial
      \G_{bc}}{\partial x^d} \right) \, .
\end{gather}
The Levi-Civita connection is torsion-free, that is $T(X,Y)=0$, for any $X$,
$Y \in \frak (M)$.

Therefore, a pseudo-Riemannian metric on $M$ defines a unique affine
connection on $M$, but the converse is not always true. Given a
pseudo-Riemannian metric $\G$ on $M$, we can define the so-called Beltrami
bracket~\cite{Cr,Sc} of functions on $M$,
\[
\Beltrami{f}{g}_\G = \G (\grad_\G f,\grad_\G g) \, , \quad f, g \in
C^\infty(M) \, .
\]
In local coordinates, one has the expression,
\[
\Beltrami{f}{g}_\G = \pder{f}{x^a} \G^{ab} \pder{g}{x^b} \, .
\]
It is interesting to note that the mapping
\[
\grad_\G : (C^{\infty}(M), \Beltrami{\cdot}{\cdot}_\G ) \rightarrow
(\mathfrak{X} (M),\symprod{\cdot}{\cdot}_{\nabla^\G})
\]
is a homomorphism of symmetric algebras, i.e., $\grad_\G \Beltrami{f}{g}_\G =
\symprod{\grad_\G f}{\grad_\G g}_{\nabla^\G}$, for all $f,g \in
C^{\infty}(M)$.

\begin{remark}
  {\rm The latter observation is the gradient analog of the following fact in
    the Hamiltonian setting: consider the mapping $(C^{\infty}(M), \{\cdot ,
    \cdot \} ) \rightarrow (\mathfrak{X} (M), [\cdot, \cdot])$ (where $\{\cdot
    , \cdot \}$ denotes the Poisson bracket and $[\cdot, \cdot]$ denotes the
    Lie bracket), associating to each function $f$ its Hamiltonian vector
    field $X_f$. Then this mapping is a homomorphism of Lie algebras, i.e.
    $X_{\{f,g\}} = [X_f,X_g]$.}
\end{remark}

Let us now turn our discussion to the cotangent bundle of $M$. First, we
introduce the construction that associates to each vector field $X$ on $M$ a
function $V^X$ on $T^*M$, defined by $V^X (x,p) = \langle p , X(x) \rangle$.
The notion of \emph{vertical lift} of a function $V$ on $M$ to a function
$\supscr{V}{v} $ on $T^*M$ is given by $\supscr{V}{v} = V \circ \pi_M$, where
$\pi_M$ is the cotangent bundle projection. An object which will play a key
role in the subsequent discussion is the \emph{Riemannian
  extension}~\cite{PaWa,YaIs} of a torsion-free affine connection. Let
$\nabla$ be a torsion-free affine connection on $M$. Then, $\nabla$ defines a
pseudo-Riemannian metric on $T^*M$, denoted $\supscr{\G}{c}$, as the unique
(0,2)-tensor on $T^*M$ which satisfies
\[
\supscr{\G}{c} (\supscr{X}{c},\supscr{Y}{c}) = - V^{\symprod{X}{Y}} \, .
\]
The fact that this single equality completely determines the Riemannian
extension $\supscr{\G}{c}$ is a consequence of the result in Proposition~4.2
in Chapter~VII of~\cite{YaIs}, which asserts that any $(0,s)$-tensor field on
$T^*M$ is univocally defined by its action on the complete lifts of vector
fields of $M$. The matrix representations of the musical isomorphisms defined
by $\supscr{\G}{c}$ in the induced local coordinates $(x^1,\ldots,
x^n,p_1,\ldots,p_n)$ on $T^*M$ are given by
\begin{align}\label{eq:local-musical}
  \flat_{\supscr{\G}{c}} \equiv
  \left(
    \begin{array}{cc}
      -2 p_c \Gamma^c_{ab} & I_n \\
      I_n & 0
    \end{array}
  \right) \, , \quad 
  \sharp_{\supscr{\G}{c}} \equiv \left(
    \begin{array}{cc}
      0 & I_n \\
      I_n & 2 p_c \Gamma^c_{ab}
    \end{array}
  \right) \, .
\end{align}
As for the gradient vector fields associated with the functions $V^X$,
$\supscr{V}{v} \in C^{\infty}(T^*M)$, $X \in \mathfrak{X}(M)$, $V \in
C^{\infty}(M)$, one has the local expressions
\begin{align}\label{eq:local-gradient}
  \grad_{\supscr{\G}{c}} V^X = X^a \pder{}{x^a} + p_a \left( \pder{X^a}{x^b} +
    2 \Gamma^a_{bc} X^c \right) \pder{}{p_b} \, , \quad \grad_{\supscr{\G}{c}}
  \supscr{V}{v} = \pder{V}{x^a} \pder{}{p_a} \, .
\end{align}

\begin{definition}\label{dfn:gradient-extension}
  The gradient extension $\Sigma^e$ of a nonlinear system $\Sigma$ of the
  form~\eqref{eq:system} with respect to a torsion-free affine connection
  $\nabla$ on $M$ is given by
  \begin{equation}\label{eq:gradient-extension}
    \Sigma^e : \left\{
      \begin{array}{l}
        \displaystyle{\dot{x}_e = \grad_{\supscr{\G}{c}} V^{g_0} (x_e) +
        \sum_{j=1}^m  u_j(t) \grad_{\supscr{\G}{c}} {V^{g_j}} (x_e) +
        \sum_{j=1}^m u_j^e (t) \grad_{\supscr{\G}{c}} \supscr{V_j}{v} (x_e)}
        \, , \\  
        y_j = \supscr{V_j}{v} (x_e) \, , \quad y_j^a = V^{g_j} (x_e) \, ,
        \quad j=1,\ldots,m \, ,  
      \end{array}
    \right.
  \end{equation}
  where $x_e = (x,p) \in T^*M$, $x_e(0) = (x_0,p_0)$, $u = (u_1,\ldots,u_m)
  \in U \subset \real^m$, and $u^e = (u_1^e,\ldots,u_m^e) \in \real^m$.
\end{definition}

\begin{remark}
  {\rm Note that the gradient extension $ \Sigma^e$ is itself a gradient
    control system.}
\end{remark}

\section{Observability of the prolongation and the gradient
  extension}\label{se:observability}

In this section, we investigate the observability properties of the prolonged
system and the gradient extension of a nonlinear system. We start by briefly
reviewing some notions such as distinguishable points and local
observability.

Let $\Yc$ denote the space of absolutely continuous functions on $M$ with
values in $\real^m$. For a nonlinear system of the form~\eqref{eq:system}, the
\emph{input-output map} $\Rc_\Sigma : M \times \Uc \rightarrow {\cal Y}$,
$\Rc_\Sigma (x_0,u(t)) = y(t)$ is defined by assigning to each initial
condition $x_0 \in M$ and any admissible control $u(t) \in \Uc$ the output of
the system,
\[
y(t) = (V_1 (x(t,x_0,u(t))),\ldots,V_m (x(t,x_0,u(t)))) \, ,
\]
where $x(t,x_0,u(t))$ denotes the solution of $\dot{x}=g_0 (x) + \sum_{j=1}^m
u_j(t) g_j(x)$ starting at $x_0$. Now, two points $x_1$, $x_2 \in M$ are said
to be \emph{indistinguishable}, $x_1 \sim x_2$, if $\Rc_\Sigma (x_1,u(\cdot))
= \Rc_\Sigma (x_2,u(\cdot))$ for any $u(\cdot) \in {\cal U}$.

\begin{definition}\label{dfn:observability}
  A system $\Sigma$ is \emph{observable} if for any $x_1$, $x_2 \in M$, one
  has that $x_1 \sim x_2 \Rightarrow x_1 = x_2$. Alternatively, for any $x_1
  \not =x_2$, there exists an admissible control such that the output
  functions resulting from the initial conditions $x(0)=x_1$, resp.
  $x(0)=x_2$, are different.  The system is \emph{locally observable} at $x_0$
  if there exists a neighborhood $\Nc$ of $x_0$ such that this holds for
  points in $\Nc$.
\end{definition}

Denote by $\Hc$ the $\real$-linear space in $C^{\infty}(M)$ spanned by the
functions of the form $\Lc_{X_1} \Lc_{X_2} \ldots \Lc_{X_s} V_j$, with $\{ X_r
\}_{r=1}^s \subset \setdef{g_i}{i=0,1,\ldots,m}$, and $j \in \{1,\ldots,m\}$.
Alternatively, we may take $X_r$ to be arbitrary elements of the
accessibility algebra corresponding to the vector fields $g_0,g_1,\ldots,g_m$.
$\Hc$ is called the \emph{observation space} of $\Sigma$. It follows from the
analyticity assumption that the system is observable if and only if $\Hc$
distinguishes points in $M$, i.e. for every $x_1$, $x_2 \in M$ with $x_1 \not
= x_2$, there exists $V \in \Hc$ such that $V(x_1) \not = V(x_2)$,
cf.~\cite{HeKr}.

\begin{proposition}[\cite{CrSc}]\label{prop:observability-prolongation}
  Consider a nonlinear system $\Sigma$ of the form~\eqref{eq:system}, with
  observation space $\Hc$.  Then, the observation space $\Hc^p$ of the
  prolongation $\Sigma^p$ is given by $\Hc^p=\supscr{\Hc}{c} +
  \supscr{\Hc}{v} $, where $\supscr{\Hc}{c} = \setdef{\supscr{V}{c}}{V \in
    \Hc}$ and $\supscr{\Hc}{v} = \setdef{\supscr{V}{v}}{V \in \Hc}$.
\end{proposition}

The following corollary is a modified statement of Corollary~3.3
in~\cite{CrSc}.
\begin{corollary}
  Assume the codistribution $d \Hc$ is of constant rank. Then the system
  $\Sigma$ is (locally) observable if and only if its prolongation is
  (locally) observable.
\end{corollary}
\begin{proof}
  Following~\cite{HeKr}, $\Sigma$ is locally observable if and only $\rk (d
  \Hc) = \dim M$. In addition, the codistribution $d \Hc$ on $M$ has constant
  rank if and only if the codistribution $d \Hc^p$ on $TM$ has constant rank.
  Therefore, $\rk (d \Hc) = \dim M$ if and only if $\rk (d \Hc^p) = \dim TM$
  if and only if $\Sigma^p$ is locally observable.  The statement regarding
  observability is proved as in Corollary~3.3 in~\cite{CrSc}.
\quad
\end{proof}

Let us turn our attention to the observability properties of the gradient
extension of a nonlinear system of the form~\eqref{eq:system}. The following
lemma will be most helpful.
\begin{lemma}\label{le:accessory}
  Let $\nabla$ be a torsion-free affine connection on a manifold $M$, and let
  $\supscr{\G}{c}$ denote its Riemannian extension to $T^*M$. Then, for any
  vector fields $X$, $Y \in \frak (M)$, and any functions $f$, $g \in
  C^{\infty}(M)$, the following identities hold
  \vspace*{.1cm}
  \begin{enumerate}
  \item $(\grad_{\supscr{\G}{c}} V^{X}) (V^Y) =
    \Beltrami{V^X}{V^Y}_{\supscr{\G}{c}} = V^{\symprod{X}{Y}} = -
    \supscr{\G}{c} (\supscr{X}{c}, \supscr{Y}{c})$.
    \vspace*{.05cm}
  \item $(\grad_{\supscr{\G}{c}} V^X ) (\supscr{f}{v} ) =
    (\grad_{\supscr{\G}{c}} \supscr{f}{v} ) (V^X) =
    \Beltrami{V^X}{\supscr{f}{v}}_{\supscr{\G}{c}} = \supscr{X(f)}{v}$.
    \vspace*{.05cm}
  \item $(\grad_{\supscr{\G}{c}} \supscr{f}{v} ) (\supscr{g}{v}) =
    \Beltrami{\supscr{f}{v}}{\supscr{g}{v}}_{\supscr{\G}{c}}=0$.
  \end{enumerate}
  \vspace*{.05cm}
\end{lemma}

\begin{proof}
  The first equality in (i) is the definition of the Beltrami bracket
  associated with $\supscr{\G}{c}$. For the second one, we resort to the
  local expressions in~\eqref{eq:local-musical} to compute
  \begin{align*}
    \Beltrami{V^X}{V^Y}_{\supscr{\G}{c}} & = \left( p_a \pder{X^a}{x^b} , X^b
    \right) \left(
      \begin{array}{cc}
        0 & I \\
        I & 2 p_e \Gamma^e_{cd}
      \end{array}
    \right) \left( p_a \pder{Y^a}{x^b} , Y^b \right)^T
    \\
    & = p_a \left( \pder{X^a}{x^b} Y^b + \pder{Y^a}{x^b} X^b + 2
      \Gamma^a_{bc} X^b Y^c \right) = V^{\symprod{X}{Y}} \, .
  \end{align*}
  The third equality corresponds to the definition of $\supscr{\G}{c}$. The
  first and second equalities in (ii) follow again by definition. As for the
  third one, note that
  \[
  \grad_{\supscr{\G}{c}} \supscr{f}{v} (V^X) = \pder{f}{x^a} \pder{}{p_a}
  (p_b X^b) = \pder{f}{x^a} X^a = \supscr{X(f)}{v}
  \]
  Finally, the equalities in (iii) are straightforward.
\quad
\end{proof}

Denote by $S_0$ the $\real$-linear space in $\frak (M)$ spanned by the vector
fields of the form
$\symprod{X_1}{\symprod{X_2}{\symprod{\ldots}{\symprod{X_{s}}{g_j}} \ldots
    }}$, with $\{ X_r \}_{r=1}^s \subset \setdef{g_i}{i=0,1,\ldots,m}$, and
$j \in \{1,\ldots,m\}$.  Alternatively, one can define $S_0$ as the smallest
subspace of $\frak (M)$ such that (i) $g_1,\ldots,g_m \in S_0$; and (ii) if
$X \in S_0$, then $\symprod{g_i}{X} \in S_0$ for all $i=0,1,\ldots,m$. We
denote by $\Sc_0$ the distribution on $M$ generated by the space $S_0$,
\[
\Sc_0 (x) = \spn \setdef{X(x)}{X \in S_0} \, , \quad x \in M \, .
\]

\begin{proposition}\label{prop:observation-extension}
  Consider a nonlinear system $\Sigma$ of the form~\eqref{eq:system}, with
  observation space $\Hc$.  Let $\nabla$ be a torsion-free affine connection
  on $M$. Then, the observation space $\Hc^e$ of the gradient extension
  $\Sigma^e$ is given by $\Hc^e=V^{\Sc_0} + \supscr{(\Hc + \mathfrak{h})}{v}$,
  where $V^{\Sc_0} = \{V^X \; | \; X \in \Sc_0\}$ and $\mathfrak{h}$ is
  spanned by $\Lc_{X_1} \Lc_{X_2} \ldots \Lc_{X_s} \Lc_X V_j$, with $X_r$,
  $r=1,\ldots,s$, equal to $g_i$, $i=0,1,\ldots,m$, $X \in \Sc_0$ and
  $j=1,\ldots,m$.
\end{proposition}

\begin{proof}
  The observation space of the gradient extension of $\Sigma$ is spanned by
  \[
  \Lc_{X_1} \Lc_{X_2} \ldots \Lc_{X_s} \supscr{V_j}{v} \, , \quad \Lc_{X_1}
  \Lc_{X_2} \ldots \Lc_{X_s} V^{g_j} \, ,
  \]
  where $X_r$, $r=1,\ldots,s$ is equal to $\grad_{\supscr{\G}{c}} V^{g_i}$,
  $\grad_{\supscr{\G}{c}} \supscr{V_j}{v}$, $i=0,1,\ldots,m$, $j=1,\ldots,m$.
  Now, using Lemma~\ref{le:accessory}, we have that
  \[
  \begin{array}{ll}
    \Lc_{\grad_{\supscr{\G}{c}} V^{g_i}} \supscr{V_j}{v} =\supscr{(\Lc_{g_i}
    V_j)}{v} \, , & \quad 
    \Lc_{\grad_{\supscr{\G}{c}} V^{g_i}} V^{g_j} = V^{\symprod{g_i}{g_j}} \,
    , \\  
    \Lc_{\grad_{\supscr{\G}{c}} \supscr{V_j}{v}} \supscr{V_k}{v} = 0 \, , &
    \quad \Lc_{\grad_{\supscr{\G}{c}} \supscr{V_j}{v}} V^{g_k} =
    \supscr{(\Lc_{g_k} V_j)}{v} \, ,  
  \end{array}
  \]
  with $i=0,1,\ldots,m$ and $j,k=1,\ldots,m$.  Considering the next step of
  Lie derivatives yields
  \[
  \begin{array}{ll}
    \Lc_{\grad_{\supscr{\G}{c}} V^{g_h}} V^{\symprod{g_i}{g_j}} =
    V^{\symprod{g_h}{\symprod{g_i}{g_j}}} \, , & \quad
    \Lc_{\grad_{\supscr{\G}{c}} V^{g_h}} \supscr{(\Lc_{g_i} V_j)}{v} =
    \supscr{(\Lc_{g_h}\Lc_{g_i} V_j)}{v}  \, , \\  
    \Lc_{\grad_{\supscr{\G}{c}} \supscr{V_k}{v}} V^{\symprod{g_i}{g_j}} =
    \supscr{(\Lc_{\symprod{g_i}{g_j}} V_k)}{v} \,  , & \quad
    \Lc_{\grad_{\supscr{\G}{c}} 
    \supscr{V_k}{v}} \supscr{(\Lc_{g_i} V_j)}{v}  = 0 \, ,  
  \end{array}
  \]
  with $h=0,1,\ldots,m$.  Further iterating this process, we get to the
  desired result.  \quad
\end{proof}

\begin{corollary}\label{cor:observability-gradient-extension}
  Consider a nonlinear system $\Sigma$ of the form~\eqref{eq:system}, with
  observation space $\Hc$.  Assume the codistribution $d \Hc$ is of constant
  rank. Let $\nabla$ be a torsion-free affine connection on $M$ and further
  assume that the distribution $\Sc_0$ is full-rank.  Then, $\Sigma$ is
  (locally) observable implies that $\Sigma^e$ is (locally) observable.
\end{corollary}

\begin{proof}
  Since the codistribution $d \Hc$ has constant rank, $\Sigma$ is locally
  observable if and only if $\dim d\Hc (x) = \dim M$.  Since $\Sc_0$ is
  full-rank, it is clear that $\Sigma$ locally observable implies that
  $\Hc^e$ has constant maximal rank, and therefore $\Sigma^e$ is locally
  observable.  With respect to observability, let $(x_1,p_1)$, $(x_2,p_2) \in
  T^*M$ and assume that $V^e (x_1,p_1) = V^e (x_2,p_2)$ for all $V^e \in
  \Hc^e$. Since $\supscr{\Hc}{v} \subset \Hc^e$, this yields $V(x_1)=V(x_2)$
  for any $V \in \Hc$. So, under observability of $\Sigma$, we conclude that
  $x_1 = x_2 = x$.  Then, we have that $V^X (x,p_1) = V^X (x,p_2)$, for all
  $X \in S_0$, which finally implies that $p_1 = p_2$.  \quad
\end{proof}

\section{Externally equivalent systems}\label{se:externally-equivalent}

In this section we introduce the notion of (weakly) externally equivalent
systems, which will be instrumental in the statement of the main result in
Section~\ref{se:main}. Consider two nonlinear systems $\alpha=1,2$, of the
form,
\[
\Sigma^\alpha : \left\{
  \begin{array}{l}
    \displaystyle{\dot{x}^\alpha = g^\alpha_0 (x^\alpha) + \sum_{j=1}^m u_j
    g^\alpha_j (x^\alpha)}   \, , \quad x^\alpha \in M^\alpha \, , \\ 
    y_j = V^\alpha_j(x^\alpha) \, , \quad j=1,\ldots,m \, , \; u =
    (u_1,\ldots,u_m) \in U \subset \real^m \, . 
  \end{array}
\right.
\]
Denote by $\Hc^\alpha$, $\alpha=1,2$, the associated observation spaces. Take
a function $H^1 \in \Hc^1$, $H^1 = \Lc_{X_1} \ldots \Lc_{X_s} V^1_j$, with $X_r
= g^1_{i_r}$, $i_r \in \{ 0,1,\ldots,m \}$, $r=1,\ldots,s$ and $j \in
\{1,\ldots,m \}$.  Consider the function in $\Hc^2$ defined by $H^2 = \Lc_{Y_1}
\ldots \Lc_{Y_s} V^2_j$, with $Y_r = g^2_{i_r}$, $r=1,\ldots,s$.  Then we say
that $H^1$ and $H^2$ formally \emph{correspond} to each other. This notion is
useful to define the concept of \emph{weakly externally equivalent} systems.

\begin{definition}\label{dfn:weakly-externally-equivalent}
  The systems $\Sigma^1$ and $\Sigma^2$ are weakly externally equivalent if
  and only if for all $x^1 \in M^1$, there exists $x^2 \in M^2$ such that
  $H^1(x^1) = H^2 (x^2)$ for all corresponding $H^1 \in \Hc^1$, $H^2 \in
  \Hc^2$, and vice versa.
\end{definition}

\begin{definition}\label{dfn:externally-equivalent}
  The systems $\Sigma^1$ and $\Sigma^2$ are externally equivalent if and only
  if for all $x^1 \in M^1$, there exists $x^2 \in M^2$ such that the
  input-output maps corresponding to $x^1$ and $x^2$ coincide, i.e.
  $\Rc_{\Sigma^1}(x^1,u(\cdot)) = \Rc_{\Sigma^2}(x^2,u(\cdot))$, for all
  $u(\cdot) \in \Uc$, and vice versa.
\end{definition}

Equivalently, $\Sigma^1$ and $\Sigma^2$ are externally equivalent if and only
if their behaviors are equal.  Clearly, if two systems are externally
equivalent, then they are weakly externally equivalent.

\begin{proposition}\label{prop:existence-unique-diffeomorphism}
  Assume that $\Sigma^1$ and $\Sigma^2$ are weakly externally equivalent,
  observable and that the codistributions $d\Hc^\alpha$, $\alpha=1,2$, have
  constant rank. Then there exists a unique diffeomorphism $\varphi: M^1
  \rightarrow M^2$ with $\varphi^* (\Hc^2) = \Hc^1$.
\end{proposition}

\begin{proof}
  Let $x^1 \in M^1$. By definition, there exists $x^2 \in M^2$ such that $H^1
  (x^1) = H^2 (x^2)$ for all corresponding $H^1 \in \Hc^1$, $H^2 \in \Hc^2$.
  Since $\Hc^2$ distinguishes points in $M^2$, it follows that $x^2$ is
  unique. Define $\varphi: M^1 \rightarrow M^2$, $\varphi (x^1) = x^2$. Using
  $\dim d\Hc^2 = \dim M^2$ and the inverse function theorem, it follows that
  $\varphi$ is smooth.  Indeed, for each $x^2 \in M^2$, there exists a
  neighborhood $V$ of $M^2$ at $x_2$ and $\dim M^2$ independent functions
  $H_1^2,\ldots,H^2_{\dim M^2}$ on $V$, such that $\varphi$ is given by
  \[
  x^2 = (H_1^2,\ldots,H^2_{\dim M^2})^{-1} (H_1^1,\ldots,H^1_{\dim M^2}) (x^1)
  \, .
  \]
  Analogously, we can construct the inverse mapping $\varphi^{-1}:M^2
  \rightarrow M^1$, making use of the fact that $\Sigma_1$ is observable,
  which concludes the proof.
\quad
\end{proof}

\begin{corollary}\label{cor:equivalent}
  Let the systems $\Sigma^1$ and $\Sigma^2$ be observable and the
  codistributions $d\Hc^\alpha$, $\alpha=1,2$, have constant rank.  Then
  $\Sigma^1$ and $\Sigma^2$ are weakly externally equivalent if and only if
  they are externally equivalent.
\end{corollary}

\begin{proof}
  We already know that if the systems are externally equivalent, then they
  are weakly externally equivalent. Conversely, assume that $\Sigma^1$ and
  $\Sigma^2$ are weakly externally equivalent.  From
  Proposition~\ref{prop:existence-unique-diffeomorphism}, we have that there
  exists a diffeomorphism $\varphi: M^1 \rightarrow M^2$ with $\varphi^*
  (\Hc^2) = \Hc^1$. Using this latter fact, and since the vector fields
  $g_0^i$, $g_j^i$ are determined by their action as derivations on
  $\Hc^\alpha$, $\alpha=1,2$, we conclude that $\varphi_* g_0^1 = g_0^2$,
  $\varphi_* g_j^1 = g_j^2$, $j=1,\ldots,m$.  \quad
\end{proof}

\begin{remark}
  {\rm The map $\varphi$ in the previous proof is called a \emph{state space
      diffeomorphism}.  }
\end{remark}

\section{Gradient realization of a nonlinear control system}\label{se:main}

This section contains the main result of the paper.  Under certain technical
conditions, Theorem~\ref{th:characterization} below characterizes when a
nonlinear control systems admits a gradient realization. Before stating this
result, we need to introduce the novel notion of \emph{compatibility} between
a nonlinear system and an affine connection.

\begin{definition}\label{dfn:compatibility}
  Let $\nabla$ be an affine connection on $M$. A nonlinear control system
  $\Sigma$ of the form~\eqref{eq:system} is compatible with $\nabla$ if and
  only if the following two conditions hold:
  \begin{description}    
  \item (a) For all vector fields $X_1,\ldots,X_{s_1}$, $Y_1,\ldots,Y_{s_2}
    \in \{g_0,g_1,\ldots,g_m\}$, and all indexes $j, k =1,\ldots,m$,
    \begin{multline*}
      \Lc_{\symprod{X_1}{\symprod{X_2}{\symprod{\ldots
              }{\symprod{X_{s_1}}{g_j}}\ldots}}}
      \left[ \Lc_{Y_1} \Lc_{Y_2}  \ldots \Lc_{Y_{s_2}} V_k \right] \\
      = \Lc_{\symprod{Y_1}{\symprod{Y_2}{\symprod{\ldots
              }{\symprod{Y_{s_2}}{g_k}}\ldots}}} \left[ \Lc_{X_1} \Lc_{X_2}
        \ldots \Lc_{X_{s_1}} V_j \right] \, .
    \end{multline*}    
  \item (b) For all vector fields $X_1,\ldots,X_{s_1}$, $Y_1,\ldots,Y_{s_2}$,
    $Z_1,\ldots,Z_{s_3} \in \{g_0,g_1,\ldots,g_m\}$, and all indexes $j, k, l
    =1,\ldots,m$,
    \begin{multline*}
      \Lc_{ \symprod{\symprod{X_1}{\symprod{X_2}{\symprod{\ldots
                }{\symprod{X_{s_1}}{g_j}}\ldots}}}{
          \symprod{Y_1}{\symprod{Y_2}{\symprod{\ldots
                }{\symprod{Y_{s_2}}{g_k}}\ldots}}}} \left[ \Lc_{Z_1} \Lc_{Z_2}
        \ldots \Lc_{Z_{s_3}} V_l \right]
      \\
      = \Lc_{\symprod{Z_1}{\symprod{Z_2}{\symprod{\ldots
              }{\symprod{Z_{s_3}}{g_l}}\ldots}}} \left[
        \Lc_{\symprod{X_1}{\symprod{X_2}{\symprod{\ldots
                }{\symprod{X_{s_1}}{g_j}}\ldots}}} \left[ \Lc_{Y_1} \Lc_{Y_2}
          \ldots \Lc_{Y_{s_2}} V_k \right] \right] \, .
    \end{multline*}
  \end{description}
\end{definition}

\begin{remark}{\rm
    Note that a locally gradient control system of the
    form~\eqref{eq:gradient-system} is compatible in the above sense with the
    Levi-Civita connection associated with the pseudo-Riemannian metric $\G$.
    Indeed, let $\symprod{\cdot}{\cdot}$, $\Beltrami{\cdot}{\cdot}$ denote,
    respectively, the symmetric product and the Beltrami bracket corresponding
    to $\nabla^\G$ and $\G$. Take $X_{r_1} = \grad V_{\alpha_{r_1}}$, $Y_{r_2}
    = \grad V_{\beta_{r_2}}$, $Z_{r_3} = \grad V_{\gamma_{r_3}}$, $r_i \in
    \{1,\ldots,s_i \}$ (which can always be written at least locally), then
    \begin{multline*}
      \Lc_{\symprod{X_1}{\symprod{X_2}{\symprod{\ldots
              }{\symprod{X_{s_1}}{g_j}}\ldots}}} \left[ \Lc_{Y_1} \Lc_{Y_2}
        \ldots \Lc_{Y_{s_2}} V_k \right] \\ = \Lc_{\grad
        \Beltrami{V_{\alpha_1}}{\Beltrami{V_{\alpha_2}}{\Beltrami{\ldots
              }{\Beltrami{V_{\alpha_{s_1}}}{V_j}}\ldots}}} \left[
        \Beltrami{V_{\beta_1}}{\Beltrami{V_{\beta_2}}{\Beltrami{\ldots
              }{\Beltrami{V_{\beta_{s_2}}}{V_k}}\ldots}}\right] \\ =
      \Lc_{\grad
        \Beltrami{V_{\beta_1}}{\Beltrami{V_{\beta_2}}{\Beltrami{\ldots
              }{\Beltrami{V_{\beta_{s_2}}}{V_k}}\ldots}}} \left[
        \Beltrami{V_{\alpha_1}}{\Beltrami{V_{\alpha_2}}{\Beltrami{\ldots
              }{\Beltrami{V_{\alpha_{s_1}}}{V_j}}\ldots}} \right]
      \\
      = \Lc_{\symprod{Y_1}{\symprod{Y_2}{\symprod{\ldots
              }{\symprod{Y_{s_2}}{g_k}}\ldots}}} \left[ \Lc_{X_1} \Lc_{X_2}
        \ldots \Lc_{X_{s_1}} V_j \right] \, ,
    \end{multline*}
    and
    \begin{multline*}
      \Lc_{\symprod{ \symprod{X_1}{\symprod{X_2}{\symprod{\ldots
                }{\symprod{X_{s_1}}{g_j}}\ldots}}}{
          \symprod{Y_1}{\symprod{Y_2}{\symprod{\ldots
                }{\symprod{Y_{s_1}}{g_k}}\ldots}}}}
      \left[ \Lc_{Z_1} \Lc_{Z_2} \ldots \Lc_{Z_{s_3}} V_l \right] \\
      = \Lc_{ \grad \Beltrami{ \Beltrami{
            V_{\alpha_1}}{\Beltrami{V_{\alpha_2}}{\Beltrami{\ldots
                }{\Beltrami{V_{\alpha_{s_1}}}{V_j}}\ldots}}}{
          \Beltrami{V_{\beta_1}}{\Beltrami{V_{\beta_2}}{\Beltrami{\ldots
                }{\Beltrami{V_{\beta_{s_2}}}{V_k}}\ldots}}}} \\
      \left[ \Beltrami{V_{\gamma_1}}{\Beltrami{V_{\gamma_2}}{\Beltrami{\ldots
              }{\Beltrami{V_{\gamma_{s_3}}}{V_l}}\ldots}} \right]
      \\
      = \Lc_{ \grad \Beltrami{
          V_{\gamma_1}}{\Beltrami{V_{\gamma_2}}{\Beltrami{\ldots
              }{\Beltrami{V_{\gamma_{s_3}}}{V_l}}\ldots}} } \\
      \left[ \Beltrami{ \Beltrami{
            V_{\alpha_1}}{\Beltrami{V_{\alpha_2}}{\Beltrami{\ldots
                }{\Beltrami{V_{\alpha_{s_1}}}{V_j}}\ldots}}}{
          \Beltrami{V_{\beta_1}}{\Beltrami{V_{\beta_2}}{\Beltrami{\ldots
                }{\Beltrami{V_{\beta_{s_2}}}{V_k}}\ldots}}} \right]
      \\
      = \Lc_{\symprod{Z_1}{\symprod{Z_2}{\symprod{\ldots
              }{\symprod{Z_{s_3}}{g_l}}\ldots}}} \left[
        \Lc_{\symprod{X_1}{\symprod{X_2}{\symprod{\ldots
                }{\symprod{X_{s_1}}{g_j}}\ldots}}} \left[ \Lc_{Y_1} \Lc_{Y_2}
          \ldots \Lc_{Y_{s_2}} V_k \right] \right] \, .
    \end{multline*}
    }
\end{remark}

\begin{remark}{\rm
    In case the distribution $\Sc_0$ is full-rank, note that property (b) in
    the above definition implies property (a) up to a constant on each
    connected component of $M$. To see this, one can use the symmetry of the
    symmetric product to deduce from (b) that
    \begin{multline*}
      \Lc_{\symprod{Z_1}{\symprod{Z_2}{\symprod{\ldots
              }{\symprod{Z_{s_3}}{g_l}}\ldots}}} \left[
        \Lc_{\symprod{X_1}{\symprod{X_2}{\symprod{\ldots
                }{\symprod{X_{s_1}}{g_j}}\ldots}}} \left[ \Lc_{Y_1} \Lc_{Y_2}
          \ldots \Lc_{Y_{s_2}} V_k \right] \right] \\
      = \Lc_{\symprod{Z_1}{\symprod{Z_2}{\symprod{\ldots
              }{\symprod{Z_{s_3}}{g_l}}\ldots}}} \left[
        \Lc_{\symprod{Y_1}{\symprod{Y_2}{\symprod{\ldots
                }{\symprod{Y_{s_2}}{g_k}}\ldots}}} \left[ \Lc_{X_1} \Lc_{X_2}
          \ldots \Lc_{X_{s_1}} V_j \right] \right] \, .
    \end{multline*}
    Now, one concludes the result from the full-rankness of $\Sc_0$.  Another
    interesting observation in this case is that the checkability of the
    compatibility condition can be performed taking a basis of vector fields
    in $S_0$ as we discuss later in Lemma~\ref{le:checkability}.}
\end{remark}

Now, we come to the main result of the paper.

\begin{theorem}\label{th:characterization}
  Let $\Sigma$ be a nonlinear control system of the form~\eqref{eq:system}.
  Let $\nabla$ be a torsion-free affine connection defined on the state
  manifold $M$.  Assume $\Sigma$ is observable with $\dim d\Hc$ constant,
  compatible with $\nabla$ and that the distribution $\Sc_0$ is full-rank.
  Then, $\Sigma$ is a locally gradient control system if and only if its
  prolonged system $\Sigma^p$ and its gradient extension $\Sigma^{e}$ are
  weakly externally equivalent.
\end{theorem}

\begin{proof}
  $\Rightarrow$) Consider a locally gradient control system $\Sigma$ on
  $(M,\G)$ (cf.~\eqref{eq:gradient-system}), together with its prolongation
  $\Sigma^p$ on $TM$ and its gradient extension $\Sigma^e$ on $T^*M$.  Recall
  that in the induced bundle coordinates $(x^a,v^a)$ on $TM$, $(x^a,p_a)$ on
  $T^*M$, the musical isomorphisms associated with $\G$ read $\flat_\G
  (x^a,v^a) = (x^a,\G_{ab}v^b)$ and $\sharp_\G (x^a,p_a) = (x^a,\G^{ab}p_b)$.
  We are going to show that $\flat_\G$ is actually an isomorphism between the
  prolongation and the gradient extension, i.e. we will prove that $\flat_\G
  (x_p (\cdot)) = x_e (\cdot)$ along the solutions
  of~\eqref{eq:prolonged-system} and~\eqref{eq:gradient-extension}
  respectively. This will be a consequence of the following equalities
  \begin{equation}\label{eq:equivalent-by-flat}
    \begin{array}{ll}
      (\flat_\G)_* \supscr{g_i}{c} = \grad_{\supscr{\G}{c}} V^{g_i} \circ
      \flat_\G \, ,  & \qquad V^{g_j} \circ \flat_\G = \supscr{V_j}{c} \, ,
      \\ 
      (\flat_\G)_* \supscr{g_j}{v} = \grad_{\supscr{\G}{c}} \supscr{V_j}{v}
      \circ \flat_\G \, ,  & \qquad  \supscr{V_j}{v} \circ \flat_\G =
      \supscr{V_j}{v} \, , 
    \end{array}
  \end{equation}
  for all $i=0,1,\ldots,m$, $j=1,\ldots,m$. In order to
  show~\eqref{eq:equivalent-by-flat}, we will make use of the following
  identities,
  \[
  (\flat_\G)_* \left( \pder{}{x^a} \right) = \pder{}{x^a} +
  \pder{\G_{cb}}{x^a} v^b \pder{}{p_c} \, , \quad (\flat_\G)_* \left(
    \pder{}{v^a} \right) = \G_{ab} \pder{}{p_b} \, .
  \]
  Let $g \in \frak (M)$. In local coordinates, $g=g^a {\partial}/{\partial
    x^a}$.  Using~\eqref{eq:complete-lift}, we get
  \[
  (\flat_\G)_* \left( \supscr{g}{c} \right) = g^a \pder{}{x^a} + \left\{ g^c
    \pder{\G_{ab}}{x^c} + \G_{ac} \pder{g^c}{x^b} \right\} v^b \pder{}{p_a} \,
  .
  \]
  On the other hand, we have that
  \[
  \grad_{\supscr{\G}{c}} V^{g} \circ \flat_\G = g^a \pder{}{x^a} + \left\{
    \G_{bc} \pder{g^c}{x^a} + 2 \G_{bc} \Gamma^c_{ad} g^d \right\} v^b
  \pder{}{p_a} \, .
  \]
  Now, suppose that $g$ is a locally gradient vector field.  In local
  coordinates, this means that $\G_{ac} g^c = \partial V / \partial x^a$, for
  a certain function $V$, which in turn implies that $\partial \{ \G_{ac} g^c
  \}/\partial x^b = \partial \{ \G_{bc} g^c \}/\partial x^a$, that is
  \begin{align*}
    \G_{ac} \pder{g^c}{x^b} = \pder{\G_{bc}}{x^a} g^c + \G_{bc}
    \pder{g^c}{x^a} - \pder{\G_{ac}}{x^b} g^c \, .
  \end{align*}
  Substituting into the above expression for $(\flat_\G)_* \left(
    \supscr{g}{c} \right)$,
  \begin{align*}
    (\flat_\G)_* \left( \supscr{g}{c} \right) & = g^a \pder{}{x^a} + \left\{
      g^c \pder{\G_{ab}}{x^c} + \pder{\G_{bc}}{x^a} g^c + \G_{bc}
      \pder{g^c}{x^a} -
      \pder{\G_{ac}}{x^b} g^c \right\} v^b \pder{}{p_a} \\
    & = g^a \pder{}{x^a} + \left\{ g^c \left( \pder{\G_{ab}}{x^c} +
        \pder{\G_{bc}}{x^a} - \pder{\G_{ac}}{x^b} \right) + \G_{bc}
      \pder{g^c}{x^a} \right\} v^b \pder{}{p_a} \\
    & = g^a \pder{}{x^a} + \left\{ 2 g^c \G_{bd} \Gamma_{ac}^d + \G_{bc}
      \pder{g^c}{x^a} \right\} v^b \pder{}{p_a} = \grad_{\supscr{\G}{c}} V^{g}
    \circ \flat_\G \, .
  \end{align*}
  Therefore, the first equality in~\eqref{eq:equivalent-by-flat} holds for
  every $i=0,1,\ldots,m$.  The equality $(\flat_\G)_* \supscr{g_j}{v} =
  \grad_{\supscr{\G}{c}} \supscr{V_j}{v} \circ \flat_{\G}$, $j=1,\ldots,m$,
  follows by considering~\eqref{eq:vertical-lift} and the fact that the
  vector fields $g_j$ are gradient by hypothesis,
  \begin{align*}
    (\flat_\G)_* \left( \supscr{g}{v} \right) = \G_{ab} g^b \pder{}{p_a} =
    \pder{V}{x^a} \pder{}{p_a} = \grad_{\supscr{\G}{c}} \supscr{V}{v} \circ
    \flat_\G \, .
  \end{align*}
  As for $V^{g_j} \circ \flat_\G = \supscr{V_j}{c}$, for each $v \in T_xM$, we
  compute $V^{g_j} \circ \flat_\G (v) = \G_{ab} v^b g_j^a = \partial V_j /
  \partial x^b \cdot v^b = < dV_j,v> = \supscr{V_j}{c}(v)$. The last equality
  follows trivially. Consequently, the prolongation and the gradient extension
  of a nonlinear system $\Sigma$ which is itself gradient are externally
  equivalent, in particular weakly externally equivalent systems.
  
  To prove the converse implication, we need some intermediate steps.  \quad
\end{proof}

\begin{lemma}\label{le:mapping}
  Let $\Sigma$ be a nonlinear system of the form~\eqref{eq:system}.  Under
  the hypothesis of Theorem~\ref{th:characterization}, assume that the
  prolongation $\Sigma^p$ and the gradient extension $\Sigma^e$ are weakly
  externally equivalent. Then there exists a unique diffeomorphism
  $\varphi:TM \rightarrow T^*M$ such that
  \begin{equation}\label{eq:equivalent-by-varphi}
    \begin{array}{ll}
      (\varphi)_* \supscr{g_i}{c} = \grad_{\supscr{\G}{c}} V^{g_i} \circ
      \varphi \, ,  & \qquad V^{g_j} \circ \varphi = \supscr{V_j}{c} \, ,
      \\ 
      (\varphi)_* \supscr{g_j}{v} = \grad_{\supscr{\G}{c}} \supscr{V_j}{v}
      \circ \varphi \, , &  \qquad \supscr{V_j}{v} \circ \varphi =
      \supscr{V_j}{v} \, , 
    \end{array}
  \end{equation}
  for all $i=0,1,\ldots,m$, $j=1,\ldots,m$. Moreover, $\varphi$ is a bundle
  morphism over the identity $\id_M : M \rightarrow M$, i.e. in natural
  coordinates $\varphi (x,v) = (x,\phi(x,v))$, for certain map $\phi:T_xM
  \rightarrow T_x^*M$, $x \in M$.
\end{lemma}

\begin{proof}
  By Proposition~\ref{prop:observability-prolongation} and
  Corollary~\ref{cor:observability-gradient-extension}, we have that both the
  prolongation and the gradient extension are observable systems. Since they
  are also weakly externally equivalent by assumption,
  Corollary~\ref{cor:equivalent} ensures that there exists a unique
  diffeomorphism $\varphi: TM \rightarrow T^*M$
  verifying~\eqref{eq:equivalent-by-varphi}.  Applying now
  Corollary~\ref{cor:equivalent} to $\Sigma^1 = \Sigma = \Sigma^2$, we deduce
  that there exists a unique diffeomorphism from $M$ to $M$ mapping the
  original nonlinear system to itself, namely the identity mapping. Using
  uniqueness and the fact that $\varphi$
  satisfies~\eqref{eq:equivalent-by-varphi}, it then follows that $\varphi$
  is of the form $\varphi (x,v) = (x,\phi(x,v))$, for certain map $\phi:T_xM
  \rightarrow T_x^*M$, $x \in M$.  \quad
\end{proof}

\begin{lemma}\label{le:mapping-structure}
  Under the same assumptions as in Lemma~\ref{le:mapping}, there exists a
  unique pseudo-Riemannian metric $\G$ on $M$ such that $\flat_{\G} =
  \varphi$, i.e. $\flat_{\G} (v) = \phi (x,v)$ for all $v \in T_xM$.
\end{lemma}

\begin{proof}
  It follows from $V^{g_j} \circ \varphi = \supscr{V_j}{c}$ (cf.
  equation~\eqref{eq:equivalent-by-varphi}) and the structure of the
  diffeomorphism $\varphi$ that
  \[
  < \phi (x,v), g_j(x) > = < dV_j (x), v > \, , \quad \forall v \in T_x M
  \, , \; j=1,\ldots,m \, .
  \]
  Furthermore, from $(\varphi)_* \supscr{g_i}{c} = \grad V^{g_i} \circ
  \varphi$ (see eq.~\eqref{eq:equivalent-by-varphi}), it follows that
  \[
  \Lc_{\grad V^{g_i}} V^{g_j} \circ \varphi = \Lc_{\supscr{g_i}{c}}
  \supscr{V_j}{c} \, , \quad i=0,1,\ldots,m \, , \; j=1,\ldots,m \, .
  \]
  Using now Lemma~\ref{le:accessory} (i), we get $<\phi
  (x,v),\symprod{g_i}{g_j}(x)> = <d \left( \Lc_{g_i}V_j \right)(x), v>$. In
  general for all $v \in T_x M$,
  \begin{multline}\label{eq:almost-varphi-adjoint}
    <\phi(x,v) , \symprod{X_1}{\symprod{X_2}{\symprod{X_3, 
          \ldots }{\symprod{X_{s}}{g_j}}\ldots}} (x) > =
    \\
    = < d \left(\Lc_{X_1} \Lc_{X_2} \ldots \Lc_{X_s} V_j \right)(x), v > \, ,
  \end{multline}
  with the $X_r$, $r=1,\ldots,s$ equal to some $g_i$, $i=0,1,\ldots,m$. Since
  the right-hand side of this equation is linear in $v$ and the distribution
  generated by the space $S_0$ is full-rank by hypothesis, it follows that
  for each $x \in M$ there exists a unique matrix $\G(x)$ such that $\phi
  (x,v) = \G (x) v$.  Since $\varphi$ is a diffeomorphism, $\G(x)$ is
  non-singular for every $x$ and depends smoothly on the base point. Consider
  the adjoint mapping of $\varphi$, $\varphi^T : TM \rightarrow T^*M$,
  defined by $<\varphi (v),w > = <v, \varphi^T (w)>$, $v,w \in T_x M$, $x \in
  M$. Then, $\varphi^T (x,v) = (x,\G^T (x) v)$. It follows
  from~\eqref{eq:almost-varphi-adjoint} that $\G (x)$ satisfies
  \begin{gather}\label{eq:varphi-adjoint}
    \varphi^T (\symprod{X_1}{\symprod{X_2}{\symprod{X_3, \ldots
          }{\symprod{X_{s}}{g_j}}\ldots}} (x) ) = d \left(\Lc_{X_1} \Lc_{X_2}
      \ldots \Lc_{X_s} V_j \right) (x) \, ,
  \end{gather}
  with the $X_r$ as above. Let us see now that $\G (x) = \G^T (x)$. Note that
  in local coordinates $(\varphi)_* \supscr{g_j}{v} = \grad_{\supscr{\G}{c}}
  \supscr{V_j}{v} \circ \varphi$ yields,
  \[
  \left(
    \begin{array}{cc}
      I & 0 \\
      \pder{}{x} \left( \G (x) v\right) & \G (x)
    \end{array} \right) \left(
    \begin{array}{cc}
      0 \\
      g_j (x)
    \end{array} \right) = \left(
    \begin{array}{cc}
      0 \\
      \left( \pder{V_j}{x} \right)^T (x)
    \end{array} \right) ,
  \]
  or, equivalently, $\G (x) g_j (x) = \left( \partial V_j / \partial x
  \right)^T (x)$, $j=1,\ldots,m$, which in intrinsic terms, can be written as
  $\varphi (g_j) = dV_j$. Now,
  \begin{gather*}
    < \varphi \left( \symprod{X_1}{\symprod{X_2}{\symprod{X_3, \ldots
            }{\symprod{X_{s_1}}{g_j}}\ldots}} \right),
    \symprod{Y_1}{\symprod{Y_2}{\symprod{Y_3, \ldots
          }{\symprod{Y_{s_2}}{g_k}}\ldots}} > \\
    = < \symprod{X_1}{\symprod{X_2}{\symprod{X_3, \ldots
          }{\symprod{X_{s_1}}{g_j}}\ldots}}, \varphi^T \left(
      \symprod{Y_1}{\symprod{Y_2}{\symprod{Y_3, \ldots
            }{\symprod{Y_{s_2}}{g_k}}\ldots}} \right) > .
  \end{gather*}
  Using~\eqref{eq:varphi-adjoint}, the latter is equal to
  \begin{multline*}
    < \symprod{X_1}{\symprod{X_2}{\symprod{X_3,  \ldots
          }{\symprod{X_{s_1}}{g_j}}\ldots}}, d \Lc_{Y_1} \Lc_{Y_2}
    \ldots \Lc_{Y_{s_2}} V_k > 
    \\
    = < d \Lc_{X_1} \Lc_{X_2} \ldots \Lc_{X_{s_1}} V_j ,
    \symprod{Y_1}{\symprod{Y_2}{\symprod{\ldots
          }{\symprod{Y_{s_2}}{g_k}}\ldots}} > \, , 
  \end{multline*}
  where in the last equality we have used the property (a) of the
  compatibility definition between the nonlinear system $\Sigma$ and the
  affine connection $\nabla$. Finally,
  \begin{multline*}
    < \varphi \left( \symprod{X_1}{\symprod{X_2}{\symprod{X_3, 
            \ldots }{\symprod{X_{s_1}}{g_j}}\ldots}}
    \right), \symprod{Y_1}{\symprod{Y_2}{\symprod{Y_3, 
          \ldots
          }{\symprod{Y_{s_2}}{g_k}}\ldots}} > \\
    = < \varphi^T \left( \symprod{X_1}{\symprod{X_2}{\symprod{X_3,
             \ldots
            }{\symprod{X_{s_1}}{g_j}}\ldots}} \right) ,
    \symprod{Y_1}{\symprod{Y_2}{\symprod{\ldots
          }{\symprod{Y_{s_2}}{g_k}}\ldots}} > \, ,
  \end{multline*}
  By the assumption on the full-rankness of the distribution $\Sc_0$, we
  conclude that
  \begin{multline*}
    \varphi \left( \symprod{X_1}{\symprod{X_2}{\symprod{X_3, \ldots
            }{\symprod{X_{s_1}}{g_j}}\ldots}} \right) =
    \\
    = \varphi^T \left( \symprod{X_1}{\symprod{X_2}{\symprod{X_3, \ldots
            }{\symprod{X_{s_1}}{g_j}}\ldots}} \right) \, ,
  \end{multline*}
  which in turn implies that $\varphi (x) = \varphi^T (x)$, i.e., the matrix
  $\G (x)$ is symmetric.  \quad
\end{proof}

\begin{lemma}\label{le:concluding}
  Under the same assumptions as in Lemma~\ref{le:mapping}, the torsion-free
  affine connection $\nabla$ is the Levi-Civita connection corresponding to
  the pseudo-Riemannian metric~$\G$.
\end{lemma}

\begin{proof}
  First of all, note that
  \begin{align}\label{eq:auxiliary-III}
    & < V^{\symprod{\symprod{X_1}{\symprod{X_2}{\symprod{\ldots
              }{\symprod{X_{s}}{g_j}} \ldots }}}{
        \symprod{Y_1}{\symprod{Y_2}{\symprod{\ldots }{\symprod{Y_{s_2}}{g_k}}
            \ldots }}}}, d (\Lc_{Z_1}
    \Lc_{Z_2} \ldots \Lc_{Z_{s_3}} V_l) > \nonumber \\
    & \quad = \Lc_{\symprod{\symprod{X_1}{\symprod{X_2}{\symprod{\ldots
              }{\symprod{X_{s}}{g_j}} \ldots }}}{
        \symprod{Y_1}{\symprod{Y_2}{\symprod{\ldots }{\symprod{Y_{s_2}}{g_k}}
            \ldots }}}} \left[
      \Lc_{Z_1} \Lc_{Z_2} \ldots \Lc_{Z_{s_3}} V_l \right] \nonumber \\
    & \quad = \Lc_{\symprod{Z_1}{\symprod{Z_2}{\symprod{\ldots
            }{\symprod{Z_{s_3}}{g_l}}\ldots}}} \left[
      \Lc_{\symprod{X_1}{\symprod{X_2}{\symprod{\ldots
              }{\symprod{X_{s}}{g_j}}\ldots}}} \left[ \Lc_{Y_1}
        \Lc_{Y_2} \ldots \Lc_{Y_{s_2}} V_k \right] \right] \\
    & \quad = < \supscr{\left(
        \Lc_{\symprod{X_1}{\symprod{X_2}{\symprod{\ldots
                }{\symprod{X_{s}}{g_j}}\ldots}}} \left[\Lc_{Y_1} \Lc_{Y_2}
          \ldots \Lc_{Y_{s_2}} V_k \right] \right)}{c} \circ \varphi^{-1}, d
    (\Lc_{Z_1} \Lc_{Z_2} \ldots \Lc_{Z_{s_3}} V_l) > . \nonumber
  \end{align}
  Since the observation space of the nonlinear system $\Sigma$ is generated by
  the functions of the form $\Lc_{Z_1} \Lc_{Z_2} \ldots \Lc_{Z_{s_3}} V_l$,
  and $\Sigma$ is observable by hypothesis, we conclude that
  \begin{multline*}
    V^{\symprod{\symprod{X_1}{\symprod{X_2}{\symprod{\ldots
              }{\symprod{X_{s}}{g_j}} \ldots }}}{
        \symprod{Y_1}{\symprod{Y_2}{\symprod{\ldots }{\symprod{Y_{s_2}}{g_k}}
            \ldots }}}} \circ
    \varphi = \\
    = \supscr{\left( \Lc_{\symprod{X_1}{\symprod{X_2}{\symprod{\ldots
                }{\symprod{X_{s}}{g_j}}\ldots}}} \left[\Lc_{Y_1} \Lc_{Y_2}
          \ldots \Lc_{Y_{s_2}} V_k \right] \right)}{c} \, .
  \end{multline*}
  Given the structure of the mapping $\varphi$ (cf. Lemmas~\ref{le:mapping}
  and~\ref{le:mapping-structure}), and equation~\eqref{eq:varphi-adjoint},
  this equality can be rewritten as,
  \begin{multline*}
    \G (\symprod{\symprod{X_1}{\symprod{X_2}{\symprod{\ldots
            }{\symprod{X_{s}}{g_j}} \ldots }}}{
      \symprod{Y_1}{\symprod{Y_2}{\symprod{\ldots
            }{\symprod{Y_{s_2}}{g_k}} \ldots }}} ,\cdot) \\
    = d < \varphi (\symprod{Y_1}{\symprod{Y_2}{\symprod{\ldots
          }{\symprod{Y_{s_2}}{g_k}} \ldots }} ) ,
    \symprod{X_1}{\symprod{X_2}{\symprod{\ldots
          }{\symprod{X_{s}}{g_j}}\ldots}} > \\
    = d \left( \G ( \grad_{\G} \left( \Lc_{Y_1} \Lc_{Y_2} \ldots
        \Lc_{Y_{s_2}} V_k \right), \grad_{\G} \left( \Lc_{X_1} \Lc_{X_2}
        \ldots \Lc_{X_{s}}
        V_j \right) ) \right) \\
    = d \Beltrami{\Lc_{Y_1} \Lc_{Y_2} \ldots \Lc_{Y_{s_2}} V_k}{\Lc_{X_1}
      \Lc_{X_2} \ldots \Lc_{X_{s}} V_j}_{\G} \, .
  \end{multline*}
  Since $\grad_\G \Beltrami{f}{g}_\G = \symprod{\grad_\G f}{\grad_\G
    g}_{\nabla^\G}$, we conclude
  \begin{multline*}
    \symprod{\symprod{X_1}{\symprod{X_2}{\symprod{\ldots
            }{\symprod{X_{s}}{g_j}} \ldots }}}{
      \symprod{Y_1}{\symprod{Y_2}{\symprod{\ldots
            }{\symprod{Y_{s_2}}{g_k}} \ldots }}} = \\
    = \symprod{\symprod{X_1}{\symprod{X_2}{\symprod{\ldots
            }{\symprod{X_{s}}{g_j}} \ldots }}}{
      \symprod{Y_1}{\symprod{Y_2}{\symprod{\ldots }{\symprod{Y_{s_2}}{g_k}}
          \ldots }}}_\G \, .
  \end{multline*}
  Using the fact that $\Sc_0$ is full-rank, we deduce that $\symprod{X}{Y} =
  \symprod{X}{Y}_\G$ for all $X$, $Y \in \mathfrak{X} (M)$. Finally, using the
  fact that $\nabla$ is torsion-free, we compute
  \begin{align*}
    \nabla_X Y = \frac{1}{2} \left( \symprod{X}{Y} + [X,Y] \right) =
    \frac{1}{2} \left( \symprod{X}{Y}_\G + [X,Y] \right) = \nabla^{\G}_X Y \,
    , \quad \forall X, Y \in \mathfrak{X} (M) \, ,
  \end{align*}
  which concludes the result.  \quad
\end{proof}

We are now ready to conclude the proof of Theorem~\ref{th:characterization}.

{\em Proof of Theorem~\ref{th:characterization}}. $\Leftarrow$) Assume the
prolongation $\Sigma^p$ and the gradient extension $\Sigma^e$ are weakly
externally equivalent. From
Lemmas~\ref{le:mapping},~\ref{le:mapping-structure} and~\ref{le:concluding},
we deduce the existence of a pseudo-Riemannian metric $\G$ on $M$ such that
$\nabla = \nabla^{\G}$ and the unique diffeomorphism between $TM$ and $T^*M$
relating $\Sigma^p$ and $\Sigma^e$ and
verifying~\eqref{eq:equivalent-by-varphi} is $\flat_{\G}$. From $(\flat_\G)_*
\supscr{g_j}{v} = \grad_{\supscr{\G}{c}} \supscr{V_j}{v} \circ \flat_{\G}$,
we deduce $\flat_\G (g_j) = dV_j$, and hence $\grad_{\G} V_j = g_j$,
$j=1,\ldots,m$. Finally, we show that $g_0$ is a locally gradient vector
field.  From $(\flat_{\G})_* \supscr{g_0}{c} = \grad_{\supscr{\G}{c}} V^{g_0}
\circ \flat_{\G}$ and the local expression~\eqref{eq:Christoffel} of the
Christoffel symbols of the Levi-Civita connection $\nabla^{\G}$, we deduce
that
\[
\pder{}{x^b} \left( \G_{ac} g_0^c \right) = \pder{}{x^a} \left( \G_{bc} g_0^c
\right) \, , \quad \forall a,b=1,\ldots,n \, ,
\]
which implies that the one-form $\flat_{\G} (g_0)$ is closed. \quad\endproof

\begin{remark}
  {\rm Note that, given the torsion-free affine connection $\nabla$,
    the pseudo-Riemannian metric $\G$ obtained in the proof of
    Theorem~\ref{th:characterization} is unique such that $\Sigma$ is
    locally gradient with respect to it. In
    Section~\ref{se:uniqueness-realization} below, we investigate the
    uniqueness (up to isometry) of gradient realizations with the same
    input-output behavior.}
\end{remark}

\begin{remark}
  {\rm In general, we cannot ensure that the drift vector field $g_0$
    is globally gradient, unless we impose some additional conditions
    on the topology of the state space $M$ (for instance, that the
    first Betti number of $M$ is zero). This is analogous to the
    situation in the Hamiltonian setting~\cite{CrSc}.}
\end{remark}

\begin{remark}
  {\rm One can verify that the pseudo-Riemannian metric on $TM$
    defined by $(\flat_\G)^* \supscr{\G}{c}$ corresponds to the
    complete lift to $TM$ of the original metric $\G$ on $M$
    (see~\cite{YaIs}).}
\end{remark}

\begin{remark}{\rm
    A different way to prove the same result which indeed keeps a closer
    parallelism with the proof for the Hamiltonian case~\cite{CrSc} would be
    the following.  Once one has proved Lemmas~\ref{le:mapping}
    and~\ref{le:mapping-structure}, instead of proving
    Lemma~\ref{le:concluding}, one can show that
    \begin{align}\label{eq:homo}
      (\varphi)_* \supscr{\symprod{X_1}{\symprod{X_2}{\symprod{\ldots
              }{\symprod{X_{s}}{g_j}}\ldots}}}{c} = \grad
      V^{\symprod{X_1}{\symprod{X_2}{\symprod{\ldots
              }{\symprod{X_{s}}{g_j}}\ldots}}} \circ \varphi \, ,
    \end{align}
    for any $j \in \{1,\ldots,m\}$ and $X_r \in \{ g_0,g_1,\ldots,g_m\}$,
    $r=1,\ldots,s$. This can be done by considering the following vector
    fields on $T^*M$,
    \begin{align*}
      \Zc_1 & = (\varphi)_*
      \supscr{\symprod{X_1}{\symprod{X_2}{\symprod{\ldots
              }{\symprod{X_{s}}{g_j}}\ldots}}}{c} \circ
      \varphi^{-1} \, , \\
      \Zc_2 & = \grad V^{\symprod{X_1}{\symprod{X_2}{\symprod{\ldots
              }{\symprod{X_{s}}{g_j}}\ldots}}} \, ,
    \end{align*}
    and showing that their action on the observation space $\Hc^e$ of
    $\Sigma^e$ is the same. To see this, recall from
    Proposition~\ref{prop:observation-extension} that $\Hc^e = V^{S_0} +
    \supscr{(\Hc + \mathfrak{h})}{v}$.  Consider a function of the form
    $\Lc_{X_1} \Lc_{X_2} \ldots \Lc_{X_s} V_j$, with $X_r$, $r=1,\ldots,s$,
    equal to $g_i$, $i=0,1,\ldots,m$, and $j=1,\ldots,m$.  Then,
    \begin{multline*}
      \Lc_{\Zc_1} \left[ \supscr{(\Lc_{X_1} \Lc_{X_2} \ldots \Lc_{X_s}
          V_j)}{v}
      \right] = \\
      = \left( \Lc_{\supscr{\symprod{X_1}{\symprod{X_2}{\symprod{\ldots
                  }{\symprod{X_{s}}{g_j}}\ldots}}}{c}} \left[
          \supscr{(\Lc_{X_1} \Lc_{X_2} \ldots \Lc_{X_s} V_j)}{v} \circ
          \varphi
        \right] \right) \circ \varphi^{-1} = \\
      = \supscr{\left( \Lc_{\symprod{X_1}{\symprod{X_2}{\symprod{\ldots
                  }{\symprod{X_{s}}{g_j}}\ldots}}} \left[ \Lc_{X_1} \Lc_{X_2}
            \ldots \Lc_{X_s} V_j \right]\right)}{V} \, ,
    \end{multline*}
    where we have used twice the fact that $\varphi$ is the identity mapping
    on the base manifold $M$. On the other hand,
    \begin{align*}
      \Lc_{\Zc_2} \left[\supscr{(\Lc_{X_1} \Lc_{X_2} \ldots \Lc_{X_s}
          V_j)}{v} \right] = \supscr{\left(
          \Lc_{\symprod{X_1}{\symprod{X_2}{\symprod{\ldots
                  }{\symprod{X_{s}}{g_j}}\ldots}}} \left[ \Lc_{X_1} \Lc_{X_2}
            \ldots \Lc_{X_s} V_j \right]\right)}{V} \, ,
    \end{align*}
    using property (ii) in Lemma~\ref{le:accessory}. The same argument also
    guarantees that the action of $\Zc_1$ and $\Zc_2$ is the same over the
    vertical lifts of the functions spanning $\mathfrak{h}$. Finally, let
    $\symprod{Y_1}{\symprod{Y_2}{\symprod{\ldots }{\symprod{Y_{s_2}}{g_k}}
        \ldots }} \in S_0$ and consider the corresponding function on $T^*M$,
    $V^{\symprod{Y_1}{\symprod{Y_2}{\symprod{\ldots }{\symprod{Y_{s_2}}{g_k}}
          \ldots }}}$. Then,
    \begin{multline}\label{eq:aux-I}
      \Lc_{\Zc_1} \left[ V^{\symprod{Y_1}{\symprod{Y_2}{\symprod{\ldots
                }{\symprod{Y_{s_2}}{g_k}} \ldots }}} \right] = \\
      = \left( \Lc_{\supscr{\symprod{X_1}{\symprod{X_2}{\symprod{\ldots
                  }{\symprod{X_{s}}{g_j}}\ldots}}}{c}} \left[
          V^{\symprod{Y_1}{\symprod{Y_2}{\symprod{\ldots
                  }{\symprod{Y_{s_2}}{g_k}} \ldots }}}
          \circ \varphi  \right] \right) \circ \varphi^{-1} = \\
      = \left( \Lc_{\supscr{\symprod{X_1}{\symprod{X_2}{\symprod{\ldots
                  }{\symprod{X_{s}}{g_j}}\ldots}}}{c}}
        \left[\supscr{(\Lc_{Y_1} \Lc_{Y_2} \ldots \Lc_{Y_{s_2}} V_k)}{c}
        \right]
      \right) \circ  \varphi^{-1} \\
      = \supscr{\left( \Lc_{\symprod{X_1}{\symprod{X_2}{\symprod{\ldots
                  }{\symprod{X_{s}}{g_j}}\ldots}}} \left[\Lc_{Y_1} \Lc_{Y_2}
            \ldots \Lc_{Y_{s_2}} V_k \right] \right)}{c} \circ \varphi^{-1}
      \, ,
    \end{multline}
    where we have used equation~\eqref{eq:varphi-adjoint}. In addition,
    \begin{align}\label{eq:aux-II}
      \Lc_{\Zc_2} \left[ V^{\symprod{Y_1}{\symprod{Y_2}{\symprod{\ldots
                }{\symprod{Y_{s_2}}{g_k}} \ldots }}} \right] =
      V^{\symprod{\symprod{X_1}{\symprod{X_2}{\symprod{\ldots
                }{\symprod{X_{s}}{g_j}} \ldots }}}{
          \symprod{Y_1}{\symprod{Y_2}{\symprod{\ldots
                }{\symprod{Y_{s_2}}{g_k}} \ldots }}}} \, ,
    \end{align}
    where we have used property (i) in Lemma~\ref{le:accessory}.  Now,
    equation~\eqref{eq:auxiliary-III} implies that~\eqref{eq:aux-I}
    and~\eqref{eq:aux-II} coincide. Therefore, $\Zc_1$ and $\Zc_2$ coincide
    over $\Hc^e$, and this concludes the proof of~\eqref{eq:homo}.
    
    Now, one can proceed by taking local coordinates $(x^1,\ldots,x^n)$ in $M$
    such that every coordinate function $x^i$ is of the form $\Lc_{X_1} \ldots
    \Lc_{X_s} V_j$ for a certain $j \in \{1,\ldots,m\}$ and certain vector
    fields $X_r \in \{ g_0,g_1,\ldots,g_m\}$, $r=1,\ldots,s$. It follows
    from~\eqref{eq:varphi-adjoint} that there exists $n$ independent vector
    fields $k^1,\ldots,k^n$ of the form
    $\symprod{X_1}{\symprod{X_2}{\symprod{\ldots
          }{\symprod{X_{s}}{g_j}}\ldots}}$ such that $\flat_{\G} (k^i) =
    dx^i$.  Finally, spelling out eq.~\eqref{eq:homo} for the vector fields
    $k^i$ and making use of the symmetry of $\G$, one obtains that the
    Christoffel symbols of the affine connection $\nabla$ are precisely given
    by~\eqref{eq:Christoffel}, which concludes the result. }
\end{remark}

\section{Uniqueness of the gradient
  realization}\label{se:uniqueness-realization}

In this section, we investigate the gradient analog of the following
well-known result for Hamiltonian systems: if two minimal Hamiltonian systems
have the same input-output map, then they are symplectomorphic~\cite{Ba3,Sc}.
We will see how the setting of Theorem~\ref{th:characterization} also
provides sufficient conditions under which a similar result holds for
gradient realizations.

In~\cite{Va}, P. Varaiya conjectured that if there exists a state space
diffeomorphism between two locally controllable gradient systems, then the
diffeomorphism is actually an isometry between the underlying
pseudo-Riemannian manifolds (see also~\cite{Ve}). Subsequently,
in~\cite{Ba,Ba2}, J. Basto Gon\c{c}alves produced an example of two locally
controllable and observable gradient systems living on the same state space
with state space diffeomorphism given by the identity mapping, where however
the Riemannian metrics are different; thus providing a counterexample to the
conjecture by Varaiya. For the sake of completeness, we review it in the
following.
\begin{example}[\cite{Ba,Ba2}]\label{ex:example} {\rm
    Consider two gradient systems $\Sigma^1$ and $\Sigma^2$ on $M^1 = M^2 =
    {\real}^4$ with Riemannian metrics $\G^1$ and $\G^2$ given respectively
    by
    \begin{align*}
      \G^1 (x_1,x_2,x_3,x_4) &= dx_1 \otimes dx_1 + e^{-x_4} dx_2 \otimes
      dx_2 + e^{-x_1} dx_3 \otimes dx_3 + e^{-x_3} dx_4 \otimes dx_4  \, ,\\
      \G^2 (x_1,x_2,x_3,x_4) &= dx_1 \otimes dx_1 + e^{-x_4} dx_2 \otimes
      dx_2 + (e^{-x_1}+e^{x_3}) dx_3 \otimes dx_3 \\
      & \qquad + e^{-x_3}(1+e^{2x_1}) dx_4 \otimes dx_4 - e^{x_1} (dx_3
      \otimes dx_4 + dx_4 \otimes dx_3) \, .
    \end{align*}
    Furthermore, let $\Sigma^1$ and $\Sigma^2$ have both zero drift vector
    fields and the same output functions given by
    \begin{align*}
      y_1 &= V_1(x) := x_1 \, , \quad y_2 = V_2(x) := x_2 + x_3 + x_4 \, .
    \end{align*}
    From the definition of $\G^1$ and $\G^2$, it easily follows that the
    input vector fields of both systems are the same, i.e.,
    \begin{align*}
      \grad_{\G^1} V_1 & = \grad_{\G^2} V_1 = \pder{}{x_1} \, , \\
      \grad_{\G^1} V_2 & = \grad_{\G^2} V_2 = e^{x_4} \pder{}{x_2} + e^{x_1}
      \pder{}{x_3} + e^{x_3} \pder{}{x _4}\, .
    \end{align*}
    Therefore, $\Sigma^1$ and $\Sigma^2$ are externally equivalent with state
    space diffeomorphism given by the identity mapping $\id : \real^4
    \rightarrow \real^4$. However, the metrics $\G^1$ and $\G^2$ are
    different, and hence the identity mapping is \emph{not} an isometry. It
    should also be noted that $\Sigma^1$ and $\Sigma^2$ are both controllable
    and observable.}
\end{example}

The following result shows that, under the hypotheses of
Theorem~\ref{th:characterization}, a state space diffeomorphism linking two
gradient systems is an isometry, \emph{provided} the state space
diffeomorphism is already known to respect the affine connections determined
by their respective pseudo-Riemannian metrics. A similar statement is already
contained in~\cite{Ba,Ba2}.  Here we make use of an argument given
in~\cite{CrSc}, p.~58 for the case of Hamiltonian systems.

\begin{proposition}\label{prop:uniqueness-realization}
  Let $\Sigma^1$ and $\Sigma^2$ be two gradient systems with state spaces
  $\left(M^1,\G^1\right)$ and $\left(M^2,\G^2\right)$, respectively. For
  $i=1,2$, assume that $\Sigma^i$ is observable with $\dim d\Hc^i$ constant,
  and that the distribution $\Sc_0^i$ is full-rank.  Furthermore, let
  $\Sigma^1$ and $\Sigma^2$ be externally equivalent with the corresponding
  state space diffeomorphism $\psi : M^1 \rightarrow M^2$ satisfying
  \begin{equation} \label{eq:conres}
    \psi _*(\nabla^{\G^1}_X Y) \circ \psi^{-1} = \nabla^{\G^2}_{\psi_* X \circ
    \psi^{-1}} (\psi_* Y \circ \psi^{-1}) \, , \quad \text{for all} \; X, Y
    \in \mathfrak{X}(M^1) \, .
  \end{equation}
  Then $\psi^* \G^2 = \G^1$, that is, $\psi$ is an isometry.
\end{proposition}

\begin{proof}
  By Lemmas~\ref{le:mapping} and~\ref{le:mapping-structure}, the map
  $\varphi^i = \flat_{\G^i}$ is the unique diffeomorphism
  satisfying~\eqref{eq:equivalent-by-varphi} for system $\Sigma^i$, $i=1,2$.
  It is easily checked that since $\Sigma^1$ and $\Sigma^2$ are externally
  equivalent with state space diffeomorphism $\psi$, then their prolongations
  $\Sigma^{1p}$ and $\Sigma^{2p}$ are externally equivalent with uniquely
  determined state space diffeomorphism given by $\psi_* : TM^1 \rightarrow
  TM^2$. Furthermore, it can be readily checked that the gradient extensions
  $\Sigma^{1e}$ and $\Sigma^{2e}$ are externally equivalent with state space
  diffeomorphism $\psi^* : T^*M^2 \rightarrow T^*M^1$, \emph{provided} $\psi$
  satisfies~\eqref{eq:conres}. This is because~\eqref{eq:conres} implies that
  $\psi^*$ respects the Riemannian extensions $(\supscr{\G}{c})^1$ and
  $(\supscr{\G}{c})^2$ determined, respectively, by the affine connections
  $\nabla^1$ and $\nabla^2$.  Hence, by the uniqueness of all these state
  space diffeomorphisms, we obtain the following commutative diagram
  \begin{figure}[htb!]
    \setlength{\unitlength}{0.9cm}
    \begin{center}
      \begin{picture}(2.3,2.3)(-0.7,0)
        \put(0,1.8){\makebox(0,0)[c]{$TM^1$}}
        \put(2.5,1.8){\makebox(0,0)[c]{$TM^2$}}
        \put(0,0){\makebox(0,0)[c]{$T^*M^1$}}
        \put(2.5,0){\makebox(0,0)[c]{$T^*M^2$}}
        
        \put(1.9,0){\vector(-1,0){1.3}}
        \put(0.6,1.8){\vector(1,0){1.3}}
        \put(0,1.4){\vector(0,-1){1}}
        \put(2.5,1.4){\vector(0,-1){1}}
        
        \put(1.3,2){\makebox(0,0)[c]{\footnotesize $\psi_*$}}
        \put(2.75,1){\makebox(0,0)[c]{\footnotesize $\varphi^2$}}
        \put(-0.25,1){\makebox(0,0)[c]{\footnotesize $\varphi^1$}}
        \put(1.3,-0.2){\makebox(0,0)[c]{\footnotesize $\psi^*$}}
        
      \end{picture}
    \end{center}
  \end{figure}
  
  \noindent that is,
  \begin{equation} \label{eq:iso}
    \psi^* \circ \varphi^2 \circ \psi_* = \varphi^1 \, .
  \end{equation}
  Recalling that $\varphi^i = \flat_{\G^i}$, $i = 1,2$, it is readily seen
  that~\eqref{eq:iso} is equivalent to
  \begin{equation} \label{eq:iso1}
    \psi^* \G^2 = \G^1 \, .
  \end{equation}
  that is, $\psi : \left(M^1,\G^1\right) \rightarrow \left(M^2,\G^2\right)$
  is an isometry. \quad
\end{proof}

\begin{remark}{\rm 
    Note that in Example~\ref{ex:example} the torsion-free connections
    determined by $\G^1$ and $\G^2$ are \emph{different}, and hence the
    identity map does not respect them.}
\end{remark}

\begin{remark}{\rm 
    Since eq.~\eqref{eq:iso} is equivalent to eq.~\eqref{eq:iso1}, one may
    also conclude that under the conditions of
    Theorem~\ref{th:characterization}, the state space diffeomorphism $\psi :
    M^1 \rightarrow M^2$ is an isometry if and only if $\psi^* : T^*M^2
    \rightarrow T^*M^1$ is a state space diffeomorphism between $\Sigma^{1e}$
    and $\Sigma^{2e}$.  }
\end{remark}

\section{Conclusions}\label{se:conclusions}

We have discussed necessary and sufficient conditions for a nonlinear control
system to be realizable as a gradient control system with respect to a
pseudo-Riemannian metric. The results rely on a suitable notion of
compatibility of the system with respect to a given affine connection, and on
the input-output behavior of the prolonged system and the gradient extension.
The symmetric product associated with an affine connection plays a key role
in the discussion. We believe that the developments in this paper do not only
give insight in the system-theoretic properties of the physically motivated
class of gradient control systems, but also shed light on the
differential-geometric properties of gradient and Lagrangian control systems.
Future work will include the investigation of equivalent characterizations in
terms of the input-output behavior of the original nonlinear system.

\section*{Acknowledgments}

The first author's work was partially supported by NSF grant CMS-0100162 and
by the European Union Training and Mobility of Researchers Program, ERB
FMRXCT-970137.

\section{Appendix}\label{se:appendix}

In this appendix we present a simplifying result concerning the compatibility
hypothesis in the statement of Theorem~\ref{th:characterization}. In general,
checking conditions (a) and (b) in the definition of compatibility between
the affine connection $\nabla$ and the nonlinear system $\Sigma$ cannot be
performed for every possible choice of vector fields in $\{g_0,g_1,\ldots,g_m
\}$ and $\{ V_1,\dots,V_m \}$. The following result shows that it is enough
to check the compatibility condition on a basis of vector fields and the
corresponding associated functions once we know that the prolongation and the
gradient extension of $\Sigma$ are weakly externally equivalent.

\begin{lemma}\label{le:checkability}
  Let $\nabla$ be a torsion-free affine connection. Assume $\Sigma$ is
  observable with $\dim d\Hc$ constant, and that the distribution $\Sc_0$ is
  full-rank. Assume the prolongation $\Sigma^p$ and the gradient extension
  $\Sigma^{e}$ of $\Sigma$ are weakly externally equivalent.  Then $\Sigma$ is
  compatible with $\nabla$ if and only if properties (a) and (b) are verified
  by a basis of vector fields in $S_0$.
\end{lemma}

\begin{proof}
  Let $R_1, \ldots, R_n$ be linearly independent vector fields of the form
  $R_i = \symprod{X^i_1}{\symprod{X^i_2}{\symprod{\ldots
        }{\symprod{X^i_{s_{i}}}{g_{j_i}}}\ldots}}$, $i=1,\ldots,n$. Let
  $V_{R_i}$ denote the function on $M$ given by $\Lc_{X_1^i} \ldots
  \Lc_{X_{s_i}^i} V_{j_i}$. From equation~\eqref{eq:varphi-adjoint}, we know
  that $\varphi^T (R_i) = d V_{R_i}$.  Assume properties (a) and (b) in the
  definition of the compatibility condition (cf.
  Definition~\ref{dfn:compatibility}) are verified by any combination of the
  vector fields $R_1,\ldots,R_n$ and the functions $V_{R_1},\ldots,V_{R_n}$.
  Let $X=\symprod{X_1}{\symprod{X_2}{\symprod{\ldots
        }{\symprod{X_{s}}{g_k}}\ldots}}$ be any element of $S_0$, and $V_X =
  \Lc_{X_1} \Lc_{X_2} \ldots \Lc_{X_{s}} V_k$ the associated function on $M$.
  Since $\Sc_0$ is full-rank, we have that $X= \sum_{i=1}^n f^i_X R_i$.  Then,
  \begin{align*}
    dV_X = d \Lc_{X_1} \Lc_{X_2} \ldots \Lc_{X_{s}} V_k = \varphi^T (X) =
    \sum_{i=1}^n f^i_X \varphi^T (R_i) = \sum_{i=1}^n f^i_X d V_{R_i} \, .
  \end{align*}
  Now, let us see that properties (a) and (b) are naturally verified by all
  possible choices of vector fields in $S_0$ and generating functions in
  $\Hc$. First,
  \begin{multline*}
    \Lc_{\symprod{X_1}{\symprod{X_2}{\symprod{\ldots
            }{\symprod{X_{s_1}}{g_j}}\ldots}}}
    \left[ \Lc_{Y_1} \Lc_{Y_2}  \ldots \Lc_{Y_{s_2}} V_k \right] \\
    = \sum_{i=1}^n f^i_Y d V_{R_i} \left( \sum_{j=1}^n f^j_X R_j \right) =
    \sum_{i,j=1}^n f^i_Y f^j_X d V_{R_j} \left( R_i \right)
    = \sum_{j=1}^n f^j_X d V_{R_j} \left( \sum_{i=1}^n f^i_X R_i \right) \\
    = \Lc_{\symprod{Y_1}{\symprod{Y_2}{\symprod{\ldots
            }{\symprod{Y_{s_2}}{g_k}}\ldots}}} \left[ \Lc_{X_1} \Lc_{X_2}
      \ldots \Lc_{X_{s_1}} V_j \right] \, ,
  \end{multline*}
  where we have used that condition (a) is verified by the vector fields
  $R_1,\ldots,R_n$ and the functions $V_{R_1},\ldots,V_{R_n}$. Secondly,
  \begin{multline}\label{eq:testing-b}
    \Lc_{\symprod{\symprod{X_1}{\symprod{X_2}{\symprod{\ldots
              }{\symprod{X_{s_1}}{g_j}}\ldots}}}{
        \symprod{Y_1}{\symprod{Y_2}{\symprod{\ldots
              }{\symprod{Y_{s_2}}{g_k}}\ldots}}}} \left[ \Lc_{Z_1} \Lc_{Z_2}
      \ldots \Lc_{Z_{s_3}} V_l \right]
    \\
    = \sum_{i=1}^n f^i_{Z} dV_{R_i} \left( \sum_{j=1}^n f^j_{\symprod{X}{Y}}
      R_j \right) = \sum_{i,j=1}^n f^i_{Z} f^j_{\symprod{X}{Y}} dV_{R_j}
    \left( R_i \right)
    \\
    = \sum_{j=1}^n f^j_{\symprod{X}{Y}} dV_{R_j} \left( \sum_{i=1}^n f^i_{Z}
      R_i \right) = < \sum_{j=1}^n f^j_{\symprod{X}{Y}} dV_{R_j}, Z> \, .
  \end{multline}
  Let us compute the coefficients $f^j_{\symprod{X}{Y}}$. We have
  \begin{multline*}
    \symprod{\symprod{X_1}{\symprod{X_2}{\symprod{\ldots
            }{\symprod{X_{s_1}}{g_j}}\ldots}}}{
      \symprod{Y_1}{\symprod{Y_2}{\symprod{\ldots
            }{\symprod{Y_{s_2}}{g_k}}\ldots}}} = \\
    \sum_{i,j=1}^n\symprod{f^i_X R_i}{f^j_Y R_j} = \sum_{i,j=1}^n \left( f^i_X
      f^j_Y \symprod{R_i}{R_j} + f^i_X R_i [f^j_Y] R_j + f^j_Y R_j
      [f^i_X] R_i \right) \\
    = \sum_{k=1}^n \left( \sum_{i,j=1}^n f^i_X f^j_Y f_{\symprod{R_i}{R_j}}^k
      + \sum_{i=1}^n f^i_X R_i [f^k_Y] + \sum_{j=1}^n f^j_Y R_j [f^k_X]
    \right) R_k \, .
  \end{multline*}
  Now, note that $\sum_{k=1}^n < f^k_{\symprod{R_i}{R_j}} dV_{R_k} , R_l > =
  \sum_{k=1}^n <f^k_{\symprod{R_i}{R_j}} dV_{R_l} , R_k >$ using condition (a)
  for the vector fields $R_1,\ldots,R_n$ and the functions
  $V_{R_1},\ldots,V_{R_n}$. Moreover, using condition (b), $\sum_{k=1}^n
  <f^k_{\symprod{R_i}{R_j}} dV_{R_l} , R_k > = < dV_{R_l} ,
  \symprod{R_i}{R_j}> = <d (dV_{R_i} [R_j]), R_l>$. Hence, $\sum_{k=1}^n
  f^k_{\symprod{R_i}{R_j}} dV_{R_k} = d (dV_{R_i} [R_j])$. On the other hand,
  \begin{multline*}
    f^i_X R_i [f^k_Y] dV_{R_k} = f^i_X <df^k_Y, R_i> dV_{R_k} \\
    = f^i_X <dV_{R_k}, R_i> df^k_Y + f^i_X \left( df^k_Y \wedge dV_{R_k}
    \right) (R_i,\cdot)
  \end{multline*}
  Since $f^i_X \left( df^k_Y \wedge dV_{R_k} \right) (R_i,\cdot) = f^i_X
  \left( d \left( f^k_Y dV_{R_k} \right) \right) (R_i,\cdot) = f^i_X \left( d
    \left( dV_{Y} \right) \right) (R_i,\cdot) = 0$, we have
  \[
  f^i_X R_i [f^k_Y] dV_{R_k} = f^i_X <dV_{R_k}, R_i> df^k_Y \, .
  \]
  Analogously, one can see that $f^j_Y R_j [f^k_X] dV_{R_k} = f^j_Y
  <dV_{R_k}, R_j> df^k_X$. Finally,
  \begin{multline*}
    \sum_{k=1}^n f^k_{\symprod{X}{Y}} dV_{R_k} = \sum_{k=1}^n \left(
      \sum_{i,j=1}^n f^i_X f^j_Y f_{\symprod{R_i}{R_j}}^k + \sum_{i=1}^n f^i_X
      R_i [f^k_Y] + \sum_{j=1}^n f^j_Y R_j [f^k_X] \right) dV_{R_k}
    \\
    = \sum_{i,j=1}^n f^i_X f^j_Y d (dV_{R_i} [R_j]) + \sum_{i,j=1}^n f^i_X
    <dV_{R_j}, R_i> df^j_Y + \sum_{i,j=1}^n f^j_Y <dV_{R_i}, R_j> df^i_X
    \\
    = d \left( \sum_{i,j=1}^n f^i_{X} f^j_{Y} dV_{R_i} [R_j] \right) = d
    (\Lc_Y \left[ V_{X} \right] ) \, .
  \end{multline*}
  Plugging this equality in~\eqref{eq:testing-b}, we get the desired result.
  \quad
\end{proof}

\end{document}